\font \fiverm=cmr5
\font \sixrm=cmr6
\font \sevenrm=cmr7
\font \eightrm=cmr8

\font \bigrm=cmr10 scaled \magstep1

\font \sevenbf=cmbx7

\font \bigbf=cmbx10 scaled \magstep1
\font \Bigbf=cmbx10 scaled \magstep2

\font \tengoth=eufm10
\font \sevengoth=eufm7
\font \fivegoth=eufm5

\newfam\gothfam
\textfont \gothfam=\tengoth
\scriptfont \gothfam=\sevengoth
\scriptscriptfont \gothfam=\fivegoth

%
%
\newfam\srmfam
\textfont \srmfam=\eightrm
\scriptfont \srmfam=\sixrm
\scriptscriptfont \srmfam=\fiverm

\font \tengothb=eufb10
\font \sevengothb=eufb7
\font \fivegothb=eufb5

\newfam\gothbfam
\textfont \gothbfam=\tengothb
\scriptfont \gothbfam=\sevengothb
\scriptscriptfont \gothbfam=\fivegothb

\font \tenmath=msbm10
\font \sevenmath=msbm7
\font \fivemath=msbm5

\newfam\mathfam
\textfont \mathfam=\tenmath
\scriptfont \mathfam=\sevenmath
\scriptscriptfont \mathfam=\fivemath
\def\math{\fam\mathfam\tenmath}
%
\parindent=0mm

%
\def\titre#1{\centerline{\Bigbf #1}\nobreak\nobreak\vglue 10mm\nobreak}

\def\paragraphe#1{\bigskip\goodbreak {\bigbf #1}\nobreak\vglue 12pt\nobreak}
\def\alinea#1{\medskip\allowbreak{\bf#1}\nobreak\vglue 9pt\nobreak}
\def\ssq{\smallskip\qquad}
\def\msq{\medskip\qquad}
\def\bsq{\bigskip\qquad}

%
%
\def\th#1{\bigskip\goodbreak {\bf Theorem #1.} \par\nobreak \sl }
\def\prop#1{\bigskip\goodbreak {\bf Proposition #1.} \par\nobreak \sl }
\def\lemme#1{\bigskip\goodbreak {\bf Lemma #1.} \par\nobreak \sl }
\def\cor#1{\bigskip\goodbreak {\bf Corollary #1.} \par\nobreak \sl }
\def\dem{\bigskip\goodbreak \it Proof. \rm}
\def\ndem{\bigskip\goodbreak \rm}
\def\qed{\par\nobreak\hfill $\bullet$ \par\goodbreak}
%
%
\def\uple#1#2{#1_1,\ldots ,{#1}_{#2}}
\def\corde#1#2-#3{{#1}_{#2},\ldots ,{#1}_{#3}}

\def\ordcorde#1#2-#3{{#1}_{#2} \le \cdots \le {#1}_{#3}}
\def\strictordcorde#1#2-#3{{#1}_{#2} < \cdots < {#1}_{#3}}
\def \restr#1{\mathstrut_{\textstyle |}\raise-6pt\hbox{$\scriptstyle #1$}}
\def \srestr#1{\mathstrut_{\scriptstyle |}\hbox to -1.5pt{}\raise-4pt\hbox{$\scriptscriptstyle #1$}}
\def \inver{^{-1}}
\def\dbar{d\!\!\hbox to 4.5pt{\hfill\vrule height 5.5pt depth -5.3pt
        width 3.5pt}}

\def\frac#1#2{{\textstyle {#1\over #2}}}

\def\R{{\math R}}
\def\C{{\math C}}

\def\N{{\math N}}
\def\Z{{\math Z}}

\def\fleche#1{\mathop{\hbox to #1 mm{\rightarrowfill}}\limits}
\def\gfleche#1{\mathop{\hbox to #1 mm{\leftarrowfill}}\limits}
\def\inj#1{\mathop{\hbox to #1 mm{$\lhook\joinrel$\rightarrowfill}}\limits}
\def\ginj#1{\mathop{\hbox to #1 mm{\leftarrowfill$\joinrel\rhook$}}\limits}
\def\surj#1{\mathop{\hbox to #1 mm{\rightarrowfill\hskip 2pt\llap{$\rightarrow$}}}\limits}
\def\gsurj#1{\mathop{\hbox to #1 mm{\rlap{$\leftarrow$}\hskip 2pt \leftarrowfill}}\limits}
\def\semi{\mathrel{\times}\kern -.85pt\joinrel\mathrel{\raise 1.07pt\hbox{${\scriptscriptstyle |}$}}}
%
%
\def \g#1{\hbox{\tengoth #1}}
\def \sg#1{\hbox{\sevengoth #1}}

\def\Cal #1{{\cal #1}}
%
%

\def \mop#1{\mathop{\hbox{\rm #1}}\nolimits}
\def \smop#1{\mathop{\hbox{\sevenrm #1}}\nolimits}

\def \mopl#1{\mathop{\hbox{\rm #1}}\limits}

%
%
\def \bib #1{\null\medskip \strut\llap{[#1]\quad}}
\def\cite#1{[#1]}
\magnification=\magstep1
\parindent=0cm
\def\titre#1{\centerline{\Bigbf #1}\vskip 16pt}
\def\paragraphe#1{\bigskip {\bigbf #1}\vskip 12pt}
\def\alinea#1{\medskip{\bf #1}\vskip 6pt}
\def\ssq{\smallskip\qquad}
\def\msq{\medskip\qquad}
\def\bsq{\bigskip\qquad}
\def\un{\hbox{\bf 1}}
\let\wt=\widetilde
\input xy
\xyoption{all}
\def\diagramme #1{\vskip 4mm \centerline {#1} \vskip 4mm}
\input epsf
\long\def\dessin#1#2#3#4{\null
                \bigskip
                \begingroup \epsfysize = #1 $$\epsfbox {#2}$$ \endgroup
                \centerline{\sevenbf Figure #3 : \sevenrm #4.}
		\bigskip
                \goodbreak}

\null
\titre{Hopf algebras,}\vskip -2mm
\titre{from basics to applications to renormalization}\vskip  -3mm
\centerline{\sl Revised and updated version, may 2006}
\bigskip
\centerline{\bigrm Dominique Manchon\footnote{$^1$}{\eightrm Universit\'e Blaise Pascal, CNRS - UMR 6620. manchon@math.univ-bpclermont.fr}}
\vskip 16 mm
\paragraphe{Introduction}
These notes are an extended version of a series of lectures given at Bogota from 2nd to 6th december 2002. They aim to present a self-contained introduction to the Hopf-algebraic techniques which appear in the work of A. Connes and D. Kreimer on renormalization in Quantum Field Theory \cite {CK1}, \cite {CK2}, \cite {BF}... Our point of view consists in revisiting a substantial part of their work in the abstract framework of connected graded Hopf algebras, i.e. Hopf algebras endowed with a compatible $\Z_+$ grading such that the degree zero component is one-dimensional.
\ssq
Chapter I contains a few elements of Hopf algebra theory which can be found in any good introductory text on the subject (\cite {Ab}, \cite {Sw}, \cite {Ka}...), as well as some basic tools from algebra which are necessary to understand the coradical filtration of a Hopf algebra. 
\ssq
Chapter II deals with connected graded and connected filtered Hopf algebras with emphasis on the convolution product. The main interest of these objects resides in the possibility to implement induction techniques with respect to the grading or the filtration~: the starting point is the particular form of the coproduct on a connected filtered Hopf algebra $\Cal H$~:
$$\Delta x=x\otimes\un+\un\otimes x+\hbox{terms involving elements of strictly
lower filtration degree}.$$
We derive the {\sl Birkhoff decomposition\/} of any linear map on $\Cal H$
with values into any commutative unital algebra $\Cal A$ which sends the unit
$\un$ to the unit $1_\Cal A$, provided the algebra $\Cal A$ is endowed with a {\sl
renormalization scheme\/}, i.e. a splitting~:
$$\Cal A=\Cal A_-\oplus\Cal A_+$$
into two subalgebras, where $\Cal A_+$ contains the unit $1_{\Cal A}$. Those
maps indeed form a group $G$ under the convolution product, and the Birkhoff
decomposition of any $\varphi\in G$ reads~:
$$\varphi=\varphi_-^{*-1}*\varphi_+,$$
where $\varphi_-$ and $\varphi_+$ are elements of $G$, $\varphi_-$ sends any
element of $\Cal H$ of positive degree to an element of $\Cal A_-$,
, and $\varphi_+$ is a map with values in
$\Cal A_+$. As an example, $\Cal A$ is the algebra of germs of meromorphic
functions at $z_0\in \C$ endowed with the {\sl minimal subtraction scheme\/},
i.e.  $\Cal A_+$ is the algebra of germs of holomorphic
functions at $z_0$, and $\Cal A_-$ is the subalgebra of polar parts $(z-z_0)\inver\C[(z-z_0)\inver]$.
\ssq
The Birkhoff decomposition respects two particular subgroups of $G$, namely
the group $G_1$ of algebra morphisms from $\Cal H$ into $\Cal A$ ({\sl
characters}), and the group $G_2$ of elements $\varphi$ of $G$ such that
$\varphi(xy)=\varphi(yx)$ ({\sl cocycles\/}). We end up this chapter with some
examples. The Birkhoff decomposition in this abstract algebraic setting has
applications beyond quantum field theory, for example in number theory where
it can be used to define multiple zeta values at nonpositive arguments
(\cite{GZ}, \cite{MP2}).
\ssq
Chapter III is devoted to Feynman graphs, which are the original example coming from quantum field theory, where D. Kreimer discovered the underlying Hopf algebra structure.
\ssq
In Chapter IV we define a bijective map $\wt R:G\to\g g$ where $\g g$ stands
for the linear space of maps on $\Cal H$ with values into $\Cal A$ which sends the unit $\un$ to $0$. It is uniquely defined by means of the equation~:
$$\varphi\circ Y=\varphi*\wt R(\varphi),$$
where $Y:\Cal H\to\Cal H$ is the natural biderivation $x\mapsto |x|.x$
associated to the graduation. We give an explicit formula for this {\sl
renormalization map\/} and for its inverse.
Finally, along the lines of \cite {CK2} we explore more closely the interplay between the renormalization map $\wt R$ and Birkhoff decomposition. When $\Cal A$ is the algebra of
germs of meromorphic functions at $z_0=0$ endowed with the minimal subtraction
scheme, the renormalization map $\wt R$ is closely related to the
$\beta$-function of \cite {CK2} as follows~: for any $\varphi\in G_1$ with
trivial positive part in the Birkhoff decomposition (i.e. $\varphi_+=u_\Cal
A\circ\varepsilon$, where $\varepsilon$ is the co-unit of $\Cal H$ and $u_\Cal
A$ is the unit map of $\Cal A$), we have~:
$$\beta(\varphi)=z\wt R(\varphi).$$
(see theorem IV.4.4).
\bsq
These notes find their origin in a workshop on renormalization in Quantum
Field Theory held in Nancy every monday during year 2000-2001. I would like to
thank all participants to this workshop, among them specially Philippe
Bonneau, Malte Henkel, Mohsen Masmoudi, Andr\'e Roux, J\'er\'emie Unterberger
and Tilmann Wurzbacher. Many thanks as well to the Mathematics Departments of
Universidad de Los Andes and Universidad La Nacional at Bogota for the
friendly atmosphere of this conference. I would like to thank Luis Fernandez
and Sylvie Paycha
for valuable remarks on previous versions of this text, as well as Kurusch
Ebrahimi-Fard, Dirk Kreimer and Li Guo for illuminating discussions.
\ssq
Enfin, un grand merci \`a toute l'\'equipe des Rencontres Math\'ematiques de
Glanon pour avoir bien voulu publier ces notes dans ce volume des
comptes-rendus.
\eject
\paragraphe{Contents~:}
\alinea{Introduction}
\alinea{I. A few elements of Hopf algebra theory}
\msq I.1. Tensor algebra
\ssq I.2. Algebras
\ssq I.3. Coalgebras
\ssq I.4. Convolution product
\ssq I.5. Bialgebras and Hopf algebras
\ssq I.6. Examples
\ssq I.7. Some properties of Hopf algebras
\alinea{II. Convolution product and regularization}
\msq II.1. Connected graded bialgebras
\ssq II.2. Connected filtered bialgebras
\ssq II.3. The convolution product
\ssq II.4. Algebra morphisms and cocycles
\ssq II.5. Birkhoff decomposition
\ssq II.6. The BCH approach to Birkhoff decomposition
\ssq II.7. Renormalized traces and characters
\ssq II.8. More on the graduation
\ssq II.9. Examples
\alinea{III. Hopf algebras of Feynman graphs}
\ssq III.1. Discarding exterior structures
\ssq III.2. Operations on Feynman graphs
\ssq III.3. The graded Hopf algebra structure
\ssq III.4. External structures
\alinea{IV. An approach to the renormalization group}
\msq III.1. The renormalization map
\ssq III.2. Inverting $\wt R$~: the scattering map
\ssq III.3. The residue
\ssq III.4. Renormalization map and Birkhoff decomposition
\medskip
\alinea{References}

\paragraphe{I. A few elements of Hopf algebra theory}
\alinea{I.1. Tensor algebra}
Let $k$ be any field, and let $A,B$ be vector spaces over $k$. The {\sl tensor product\/} $A\otimes B$ is a vector space over $k$ which satisfies the following {\sl universal property\/}~: there exists a bilinear map~:
$$\eqalign{\iota~:A\times B	&\longrightarrow A \otimes B	\cr
		(a,b)		&\longmapsto a\otimes b	\cr}$$
such that for any $k$-vector space $C$ and for any bilinear map $f$ from $A\times B$ into $C$ there is a unique linear map $\wt f~:A\otimes B\to C$ such that $f=\wt f\circ\iota$~:
\diagramme{
\xymatrix{
A\otimes B \ar [dr]^{\wt f}&  \\
A\times B \ar[u]^{\iota} \ar [r]^f	& C }
}
\prop{I.1.1}
Tensor product $A\otimes B$ exists and is unique up to isomorphism.
\dem
We show uniqueness first~: if $(T_1,\iota_1)$ and $(T_2,\iota_2)$ are both candidates for playing the role of the tensor product, then the universal property applied to both tells us that there exist $\phi~:T_1\to T_2$ and $\psi~:T_2\to T_1$ such that $\iota_2=\phi\circ\iota_1$ and $\iota_1=\psi\circ\iota_2$.
\diagramme{
\xymatrix{
T_1	\ar@<.5ex>[dr]^{\phi}	&	\\
A\times B	\ar[u]^{\iota_1} \ar[r]_{\iota_2}	&T_2 \ar@<.5ex>[ul]^\psi}
}
Applying the universal property twice again shows that $\psi\circ\phi=I_{T_1}$ and $\phi\circ\psi=I_{T_2}$, whence uniqueness of the tensor product up to linear isomorphism.
\ssq
For the existence we shall need axiom of choice in the infinite-dimensional case~: let $(e_i)_{i\in E}$ and $(f_j)_{j\in F}$ be bases of $A$ and $B$ respectively. Then we take as $A\otimes B$ the vector space freely generated by $(c_{ij})_{i\in E, j\in F}$, and we set $\iota(e_i,f_j)=e_i\otimes f_j=c_{ij}$. For any bilinear map $\phi$ from $A\times B$ to another vector space $C$ it is clear that the linear map $\wt \phi$ from $A\otimes B$ to $C$ defined by~:
$$\wt \phi(c_{ij})=\phi(e_i,f_j)$$
is the only one such that $\wt \phi\circ\iota=\phi$. The vector space $A\otimes B$ constructed above fulfills then the universal property.
\qed 
Elements $a\otimes b$ in $A\otimes B$ are called {\sl indecomposable\/}. They clearly generate $A\otimes B$. Tensor products $k\otimes A$ and $A\otimes k$ are canonically identified with $A$ via $1\otimes a\simeq a\simeq a\otimes 1$. When three vector spaces $A,B,C$ are involved there is an isomorphism~:
$$\eqalign{ \alpha~:(A\otimes B)\otimes C	
	&\longrightarrow A\otimes (B\otimes C)	\cr
(a\otimes b)\otimes c	&\longmapsto a\otimes (b\otimes c).\cr}$$
Note that this isomorphism is {\sl not\/} canonical, because tensor product $A\otimes B$ itself is only defined up to isomorphism. We shall denote by $A\otimes B\otimes C$ any of these two versions of the iterated tensor product.
\ssq
Let $A$ and $B$ be two vector spaces. The {\sl flip\/} $\tau:A\otimes B\to B\otimes A$ defined by $\tau(a\otimes b)=b\otimes a$ is an isomorphism. This generalizes to any finite collection of vector spaces as follows~: any permutation $\sigma\in S_n$ defines a map~:
$$\eqalign{\tau_\sigma:A_1\otimes\cdots\otimes A_n	&\longrightarrow 
A_{\sigma \inver 1}
	\otimes\cdots\otimes A_{\sigma\inver n}	\cr
	a_1\otimes\cdots\otimes a_n		&\longmapsto	a_{\sigma\inver 1}
	\otimes\cdots\otimes a_{\sigma\inver n}.	\cr}$$
In particular the group $S_n$ acts on the tensor power $A^{\otimes n}$ by automorphisms. It makes then sense to consider the tensor power $A^{\otimes F}$ for any finite set, by choosing any arbitrary total order on $F$.		
\smallskip
Now let us consider four $k$-vector spaces $A_1,A_2,B_1,B_2$.
\prop{I.1.2}
There is a canonical injection~:
$$\Cal L(A_1,A_2)\otimes\Cal L(B_1,B_2)\inj 6^j_\sim
\Cal L(A_1\otimes B_1,\,A_2\otimes B_2)$$
given by~:
$$\bigl(j(f\otimes g)\bigr)(a\otimes b)=f(a)\otimes g(b).$$
In the case when all the spaces are finite-dimensional, this injection is an isomorphism.
\dem
The space $\Cal L(A_1\otimes B_1,\,A_2\otimes B_2)$ together with the bilinear map~:
$$\eqalign{\iota~: \Cal L(A_1,A_2)\times\Cal L(B_1,B_2)	&\longrightarrow
\Cal L(A_1\otimes B_1,\,A_2\otimes B_2)	\cr
	(f,g)	&\longmapsto \bigl(a\otimes b\mapsto f(a)\otimes g(b)\bigr)
\cr}$$
yields the map $j$, which is easily seen to be injective. In the finite-dimensional case we can either compute the dimensions or verify that the space $\Cal L(A_1\otimes B_1,\,A_2\otimes B_2)$ together with the bilinear map $\iota$ fulfills the universal property. The details are left to the reader.
\qed
{\sl Remark\/}~: There are ``super'' versions of the tensor product for super-vector spaces (i.e. $\Z_2$-graded vector spaces). Elements of vector spaces and linear morphisms are decomposed into degree zero (even) and degree one (odd) components. Then proposition I.1.2 must be modified by a minus sign if both $g$ and $a$ are odd~: this is the {\sl Koszul rule of signs\/}~:
$$(f\otimes g)(a\otimes b)=(-1)^{|g||a|}f(a)\otimes g(b),$$
according to which any transposition in such a formula brings a sign which is minus if and only if the two elements involved in the transposition are odd. We shall not use this $\Z_2$-graded framework in the sequel.

\alinea{I.2. Algebras and modules}
{\sl I.2.1. Basic definitions\/}
\medskip
A $k$-algebra is by definition a $k$-vector space $A$ together with a bilinear map $m:A\otimes A\to A$ which is {\sl associative\/}. The associativity is expressed by the commutativity of the following diagram~:
\diagramme{
\xymatrix{
A\otimes A\otimes A \ar[d]^{I\otimes m}\ar[r]^{m\otimes I}
	&A\otimes A \ar[d]^m	\\
A\otimes A \ar[r]^m	& A}
}
The algebra $A$ is {\sl unital\/} if moreover there is a unit $\un$ in it. This is expressed by the commutativity of the following diagram~:
\diagramme{
\xymatrix{
k\otimes A \ar[r]^{u\otimes I} \ar[dr]^\sim	& A\otimes A \
 \ar[d]_m	&A\otimes k \ar[l]_{I\otimes u}\ar[dl]_\sim	\\
&A&}
}
where $u$ is the map from $k$ to $A$ defined by $u(\lambda)=\lambda\un$. The algebra $A$ is {\sl commutative\/} if $m\circ \tau=m$, where $\tau~:A\otimes A\to A\otimes A$ is the {\sl flip\/}, defined by $\tau(a\otimes b)=b\otimes a$.
\ssq
A subspace $J\subset A$ is called a {\sl subalgebra\/} (resp. a {\sl left ideal, right ideal, two-sided ideal\/}) of $A$ if $m(J\otimes J)$ (resp. $m(J\otimes A)$, $m(A\otimes J)$, $m(J\otimes A+A\otimes J$) is included in $J$.
\medskip
{\sl I.2.2. Algebras and tensor product\/}
\medskip
To any vector space $V$ we can associate its {\sl tensor algebra\/} $T(V)$. As a vector space it is defined by~:
$$T(V)=\bigoplus_{k\ge 0}V^{\otimes k},$$
with $V^{\otimes 0}=k$ and $V^{\otimes k+1}:=V\otimes V^{\otimes k}.$ The product is given by the {\sl concatenation\/}~:
$$m(v_1\otimes\cdots\otimes v_p,\, v_{p+1}\otimes\cdots\otimes v_{p+q})=
v_1\otimes\cdots\otimes v_{p+q}.$$
The embedding of $k=V^{\otimes 0}$ into $T(V)$ gives the unit map $u$. Tensor algebra $T(V)$ is also called the {\sl free (unital) algebra generated by $V$\/}. This algebra is characterized by the following universal property~: for any linear map $\varphi$ from $V$ to an algebra $A$ there is a unique algebra morphism $\wt\varphi$ from $T(V)$ to $A$ extending $\varphi$. The fact that this property characterizes $T(V)$ up to isomorphism is an easy exercise left to the reader.
\ssq
Let $A$ and $B$ be unital $k$-algebras. We put a unital algebra structure on $A\otimes B$ in the following way~:
$$(a_1\otimes b_1).(a_2\otimes b_2)=a_1a_2\otimes b_1b_2.$$
The unit element $\un_{A\otimes B}$ is given by $\un_A\otimes \un_B$, and the associativity is clear. This multiplication is given by~:
$$m_{A\otimes B}=(m_A\otimes m_B)\circ \tau_{23},$$
where $\tau_{23}~: A\otimes B\otimes A\otimes B\to A\otimes A\otimes B\otimes B$ is defined by the flip of the two middle factors~:
$$\tau_{23}(a_1\otimes b_1\otimes a_2\otimes b_2)=a_1\otimes a_2\otimes b_1\otimes b_2.$$
\medskip
{\sl I.2.3. Modules\/}
\medskip
Let $A$ be any unital algebra. A {\sl left $A$-module\/} is a $k$-vector space $M$ together with a map $\alpha:A\otimes M\to M$ such that the following diagrams commute~:
\diagramme{
\xymatrix{
A\otimes A\otimes M \ar[d]^{I\otimes \alpha}\ar[r]^{m\otimes I}
	&A\otimes M \ar[d]^\alpha	&&k\otimes M \ar[r]^{u\otimes I} \ar[dr]^\sim	& A\otimes M \
 \ar[d]^\alpha\\
A\otimes M \ar[r]^\alpha	& M &&&A}
}
The map $\alpha$ is called the action of the algebra $A$ on $M$. For any $a\in A$ and $m\in M$ we usually denote by $a.m$ the action $\alpha (a\otimes m)$ of $a$ on $m$. The two diagrams above express the identities~:
$$(a.b).m=a.(b.m),\hskip 20mm \un .m=m$$
for any $a,b\in A$ and $m\in M$. The {\sl right $A$-module\/} are defined similarly, replacing $A\otimes M$ with $M\otimes A$ (details are left to the reader). A linear subspace $N$ of a left $A$-module $M$ is called a {\sl submodule\/} if $\alpha(A\otimes N)\subset N$. The intersection of all left submodules of $M$ containing a subset $P$ is called the {\sl left submodule generated by $P$\/}.
\ssq
A left module $M$ is {\sl simple\/} if it does not contain any submodule different from $\{0\}$ or $M$ itself. If a left module $M$ can be written as a direct sum of simple modules, $M$ is said to be {\sl semi-simple\/}.
\prop{I.2.1}
For any left maximal ideal $J$ of an algebra $A$ the quotient $A/J$ is a simple left $A$-module, and conversely any simple left $A$-module is isomorphic to a
simple left $A$-module of this form.
\dem
The first assertion is immediate. Conversely, let $M$ a simple left module and let $m\in M-\{0\}$. Let $J_m$ be the annihilator of $m$. It is a left ideal of $A$, and by simplicity of $M$ the map~:
$$\eqalign{\phi_m~:A	&\longrightarrow M	\cr
		a	&\longmapsto a.m	\cr}$$
gives rise to an morphism of left $A$-modules from $A/J_m$ onto $M$. It is injective by definition of $J_m$, and surjectivity comes from the simplicity of the module $M$. So $\phi_m$ is an isomorphism. The left ideal $J_m$ is maximal, which proves the proposition.
\qed
Now let $M$ be an $A$-module. We denote by $A'_M$ the algebra of endomorphisms of $M$ as an $A$-module, and we denote by $A''_M$ the algebra of endomorphisms of $M$ as an $A'_M$-module. Clearly any $a\in A$ gives rise to an element of $A''_M$. Following N. Jacobson \cite {J Chap. 4.3} we shall give a proof of an important {\sl density theorem\/}~:
\th{I.2.2}
Let $M$ be a semi-simple $A$-module, and let $\uple xn$ a finite collection of elements of $M$. Then for any $a''\in A''_M$ there exists an element $a\in A$ such that $ax_i=a''x_i$ for any $i=1,\ldots ,n$.
\dem
First notice that any $A$-submodule $N$ of $M$ is an $A''_M$-submodule. To see this write (thanks to semi-simplicity) $M=N\oplus T$ where $T$ is another $A$-submodule of $M$. The projection $e$ on $N$ with respect to this decomposition is an element of $A'_M$. For any $a''\in A''_M$ we have then~:
$$a''(N)=a''\circ e(M)=e\circ a''(M)\subset N.$$
Consider then for any fixed positive integer $n$ the semi-simple module $M^n$, direct sum of $n$ copies of $M$. The algebra $A'_{M^n}$ coincides with the algebra of $n\times n$ matrices over $A'_M$, and the diagonal matrices over $A''_M$ form a subalgebra of $A''_{M^n}$, and thus realize an embedding of $A''_M$ into $A''_{M^n}$.
\ssq
Consider $x=(\uple xn)\in M^n$. Then $N=A.x$ is an $A$-submodule of $M^n$. So it is an $A''_{M^n}$-submodule, hence an $A''_M$-submodule via the diagonal embedding above. Then for any $a''\in A''_M$ there exists $a\in A$ such that $a''x=ax$, which proves the theorem.
\qed
\cor{I.2.3}
On a semi-simple finite-dimensional module $M$ the natural map from $A$ into $A''_M$ is surjective. 
\ndem
\medskip
{\sl I.2.4. The Jacobson radical\/}
\medskip
Let $A$ be a $k$-algebra. The {\sl radical\/} $\mop{rad} M$ of a left module is by definition the intersection of all maximal submodules of $M$. When the module $M$ is the algebra $A$ itself, the radical $\mop{rad} A$ is the intersection of all maximal left ideals, and is called the {\sl Jacobson radical\/} of the algebra $A$.
\ssq
We shall give an alternative definition of the Jacobson radical~: a {\sl primitive ideal\/} is the annihilator of a simple module. In view of proposition I.2.1, any primitive ideal is the annihilator of $A/J$ where $J$ is a maximal left ideal. Of course a primitive ideal is two-sided.
\lemme{I.2.4}
Any primitive ideal is an intersection of maximal left ideals.
\dem
Any primitive ideal $J$ is by definition the annihilator of a simple module $M$. The annihilator $J_m$ of any $m\in M-\{0\}$ is then a maximal left ideal containing $J$, and it is clear that we have~:
$$J=\bigcap_{m\in M-\{0\}} J_m.$$
\qed
\prop{I.2.5}
The Jacobson radical of $A$ is the intersection of its primitive ideals.
\dem
Let us call $P$ the intersection of all primitive ideals of $A$. By lemma I.2.4 and proposition I.2.1, $P$ is indeed the intersection of all maximal left ideals.
\qed
\lemme{I.2.6 \rm (Nakayama's lemma)}
Let $M$ a finitely generated $A$-module, and let $N,L$ two submodules of $M$ such that $M=L+N$ and $N\subset \mop{rad}M$. Then $L=M$.
\dem
Suppose that $L$ is strictly contained in $M$. As $M$ is finitely generated there exists a maximal nontrivial submodule $\wt L$ containing $L$. It contains $N$ as well by definition of $\mop {rad} M$. Then $\wt L$ contains $L+N$, so $L+N$ cannot be equal to $M$.
\qed
\cor{I.2.7}
The Jacobson radical of a finite-dimensional algebra is nilpotent.
\dem
Let $A$ a finite-dimensional algebra with radical $R$. Observe first that for any $A$-module $M$ we have the inclusion~:
$$R.M\subset\mop{rad} M.$$
Indeed any maximal submodule $N$ of $M$ contains $J.M$ where $J$ is a primitive ideal, namely the annihilator of $M/N$. Hence any maximal submodule of $M$ contains $R.M$. Suppose now that $A$ is finite dimensional and that $J$ is an ideal of $A$ such that $R.J=J$. A fortiori $\mop{rad}J=J$. Applying Nakayama's lemma I.2.6 to $M=J$ and $L=\{0\}$ we get $J=\{0\}$. We immediately deduce from this fact that for any positive integer $n$, $R^n$ either contains strictly $R^{n+1}$ or is equal to $\{0\}$. As $A$ is finite-dimensional $R^n$ is indeed equal to $\{0\}$ for some $n$.
\qed
{\sl Remark\/}~: although the definition of the Jacobson ideal is not symmetric (because we used left ideals and left modules), the Jacobson ideal itself is a symmetric notion~: in other words the Jacobson radicals of algebras $A$ and $A^{\smop{opp}}$ coincide. In order to see this one can show that $\mop{Rad} A$ is the biggest two-sided ideal $J$ such that $1-x$ is invertible for any $x\in J$ \cite {B \S 6.3}. This definition is indeed symmetric. 
\medskip
{\sl I.2.5. Maximal two-sided ideals\/}
\medskip
It is easily seen that any maximal two-sided ideal in an algebra $A$ is primitive. The converse is false in general : for example, in the enveloping algebra of the non-trivial two-dimensional Lie algebra, the ideal $\{0\}$ is primitive but not maximal \cite {Di}. However we have the following result~:
\prop{I.2.8}
Any finite-codimensional primitive ideal is maximal.
\ndem
Before giving a proof of this result, we need the following definitions~: an algebra $A$ is {\sl simple\/} if it does not contain any proper two-sided ideal. A {\sl division algebra\/} is an algebra such that any nonzero element is inversible. %
\lemme{I.2.9}
Let $D$ be a division algebra. Then the algebra $M_n(D)$ of square $n\times n$ matrices over $D$ is simple.
\dem
Let us consider for any $i,j\in\{1,\ldots ,n\}$ the elementary matrix $e_{ij}$ with vanishing entries except the one on the i$^{\hbox{\sevenrm th}}$ row and j$^{\hbox{\sevenrm th}}$ column which is equal to the unit of $D$. Denote by $I_i$ the left ideal of $M_n(D)$ consisting of matrices such that all columns vanish except the i$^{\hbox{\sevenrm th}}$ column. These left ideals are all simple and isomorphic as $M_n(D)$-modules. Now let $I$ a nonzero two-sided ideal of $M_n(D)$. Let $X$ a nonzero element of $I$, and $x_{kl}$ a nonzero entry of the matrix $X$. Then for any $i\in\{0,\ldots ,n\}$ the product $e_{ik}X$ belongs to $I_i\cap I$ and is different from $0$.
\ssq
Then $I_i\cap I\not =\{0\}$. As $I_i$ is simple as a left module for each $i$ that means that $I$ contains all the left ideals $I_i$, and hence $I=M_n(D)$.
\qed
There is a converse to this result, namely any simple {\sl artinian\/} algebra (a fortiori any simple finite-dimensional algebra) is isomorphic to $M_n(D)$ where $D$ is a division algebra. This is a particular case of the {\sl Wedderburn-Artin theorem\/}, which gives a complete description of {\sl semi-simple algebras\/} \cite{DF Chap. 1}, \cite{DK Chap. 2}.
\medskip
{\sl Proof of Proposition I.2.8\/}~: a finite-dimensional primitive ideal $I$ is the annihilator of a simple finite-dimensional module $M$. By simplicity of $M$ the algebra $D=A'_M$ is a division algebra (Schur's lemma). The action of $A$ on $M$ yields (thanks to corollary I.2.3) a surjective algebra morphism from $A$ onto $A''_M$, and hence an algebra isomorphism from $A/J$ onto $A''_M$. But $A''_M$ is a matrix algebra over $D$, and then is simple (according to lemma I.2.7). So $A/J$ is a simple algebra, which amounts to say that $J$ is maximal as a two-sided ideal.
\qed
\cor{I.2.10}
In a finite-dimensional algebra, primitive ideals and maximal two-sided ideals coincide. In particular the Jacobson radical is the intersection of all maximal two-sided ideals in this case.
\ndem
\alinea{I.3. Coalgebras and comodules}
This paragraph is mostly borrowed from M.E. Sweedler's book \cite {Sw}, particularly chapters 1, 2, 8 and 9.
\medskip
{\sl I.3.1. Coalgebras\/}
\medskip
Coalgebras are objects wich are somehow dual to algebras~: axioms for coalgebras are derived from axioms for algebras by reversing the arrows of the corresponding diagrams~:
\ssq
A $k$-coalgebra is by definition a $k$-vector space $C$ together with a bilinear map $\Delta:C\to C\otimes C$ which is {\sl co-associative\/}. The co-associativity is expressed by the commutativity of the following diagram~:
\diagramme{
\xymatrix{
C\otimes C\otimes C 
	&C\otimes C \ar[l]_{\Delta\otimes I}	\\
C\otimes C \ar[u]_{I\otimes \Delta}	& C\ar[l]_{\Delta}\ar[u]_{\Delta}}
}
Coalgebra $C$ is {\sl co-unital\/} if moreover there is a co-unit $\varepsilon$ such that the following diagram commutes~:
\diagramme{
\xymatrix{
k\otimes C  	& C\otimes C \ar[l]_{\varepsilon\otimes I}\
 	\ar[r]^{I\otimes \varepsilon} &C\otimes k 	\\
&C\ar[u]^\Delta \ar[ul]_\sim \ar[ur]^\sim &}
}
A subspace $J\subset C$ is called a {\sl subcoalgebra\/} (resp. a {\sl left
coideal, right coideal, two-sided coideal\/}) of $C$ if $\Delta(J)$ is
contained in $J\otimes J$ (resp. $C\otimes J$,  $J\otimes C$, $J\otimes C+C\otimes J$) is included in $J$. The duality alluded to above can be made more precise~:
\prop{I.3.1}
1)The linear dual $C^*$ of a co-unital coalgebra $C$ is a unital algebra, with product (resp. unit map) the transpose of the coproduct (resp. of the co-unity). 
\smallskip
2)Let $J$ a linear subspace of $C$. Denote by $J^\perp$ the orthogonal of $J$ in $C^*$. Then~:
\smallskip
\quad
$J$ is a two-sided coideal if and only if $J^\perp$ is a subalgebra of $C^*$.

\quad
$J$ is a left coideal if and only if $J^\perp$ is a left ideal of $C^*$.

\quad
$J$ is a right coideal if and only if $J^\perp$ is a right ideal of $C^*$.

\quad
$J$ is a subcoalgebra if and only if $J^\perp$ is a two-sided ideal of $C^*$.
\dem
For any subspace $K$ of $C^*$ we shall denote by $K^\perp$ the subspace of those elements of $C$ on which any element of $K$ vanishes. It coincides with the intersection of the orthogonal of $K$ with $C$, via the canonical embedding $C\inj 6 C^{**}$. So we have for any linear subspaces $J\subset C$ and $K\subset C^*$~:
$$J^{\perp\perp}=J,\hskip 20mm K^{\perp\perp}\supset K.$$  
Suppose that $J$ is a two-sided coideal. Take any $\xi,\eta$ in $J^\perp$. For any $x\in J$ we have~:
$$<\xi\eta,x>=<\xi\otimes\eta,\Delta x>=0,$$
as $\Delta x\subset J\otimes C+C\otimes J$. So $J^\perp$ is a subalgebra of $C^*$. Conversely if $J^\perp$ is a subalgebra then~:
$$\Delta J\subset (J^\perp\otimes J^\perp)^\perp=J\otimes C+C\otimes J,$$
which proves the first assertion. We leave it to the reader as an exercice to prove the three other assertions along the same lines. Dually we have the following~:
\prop{I.3.2}
Let $K$ a linear subspace of $C^*$. Then~:
\smallskip
\quad
$K^\perp$ is a two-sided coideal if and only if $K$ is a subalgebra of $C^*$.

\quad
$K^\perp$ is a left coideal if and only if $K$ is a left ideal of $C^*$.

\quad
$K^\perp$ is a right coideal if and only if $K$ is a rightt ideal of $C^*$.

\quad
$K^\perp$ is a subcoalgebra if and only if $K$ is a two-sided ideal of $C^*$.
\ndem
The linear dual $(C\otimes C)^*$ naturally contains the tensor product $C^*\otimes C^*$. Take as a multiplication the restriction of $^t\!\Delta$ to $C^*\otimes C^*$~:
$$m=^t\!\Delta~:C^*\otimes C^*\longrightarrow C^*,$$
and put $u=^t\!\varepsilon~:k\to C^*$. It is easily seen, by just reverting the arrows of the corresponding diagrams, that coassociativity of $\Delta$ implies associativity of $m$, and that the co-unit property for $\varepsilon$ implies that $u$ is a unit.
\qed
The coalgebra $C$ is {\sl cocommutative\/} if $\tau\circ\Delta=\Delta$, where $\tau~:C\otimes C\to C\otimes C$ is the flip. It will be convenient to use {\sl Sweedler's notation\/}~:
$$\Delta x=\sum_{(x)}x_1\otimes x_2.$$
Cocommutativity expresses then as~:
$$\sum_{(x)}x_1\otimes x_2=\sum_{(x)}x_2\otimes x_1.$$
Coassociativity reads in Sweedler's notation~:
$$(\Delta\otimes I)\circ\Delta(x)=\sum_{(x)}x_{1:1}\otimes x_{1:2}\otimes x_2=
	\sum_{(x)}x_1\otimes x_{2:1}\otimes x_{2:2}
	=(I\otimes\Delta)\circ\Delta(x),$$
We shall sometimes write the iterated coproduct as~:
$$\sum_{(x)}x_1\otimes x_2\otimes x_3.$$
Sometimes we shall even mix the two ways of using Sweedler's notation for the iterated coproduct, in the case we want to keep partially track of how we have constructed it \cite{DNR}. For example,
$$\eqalign{\Delta_3(x)	&=(I\otimes\Delta\otimes I)\circ(\Delta\otimes I)
	\circ\Delta(x)	\cr
			&=(I\otimes\Delta\otimes I)(\sum_{(x)}
		x_1\otimes x_2\otimes x_3)	\cr
			&=\sum_{(x)}
		x_1\otimes x_{2:1}\otimes x_{2:2}\otimes x_3.	\cr}$$	
\ssq
To any vector space $V$ we can associate its {\sl tensor coalgebra\/} $T^c(V)$. It is isomorphic to $T(V)$ as a vector space. The coproduct is given by the {\sl deconcatenation\/}~:
$$\Delta(v_1\otimes\cdots\otimes v_n)=\sum_{p=0}^n
(v_1\otimes\cdots\otimes v_p) \bigotimes (v_{p+1}\otimes\cdots\otimes v_n).$$
The co-unit is given by the natural projection of $T^c(V)$ onto $k$.
\ssq
Let $C$ and $D$ be unital $k$-coalgebras. We put a co-unital coalgebra structure on $C\otimes D$ in the following way~: the comultiplication is given by~:
$$\Delta_{C\otimes D}=\tau_{23}\circ(\Delta_C\otimes \Delta_D) ,$$
where $\tau_{23}$ is again the flip of the two middle factors, and the co-unity is given by $\varepsilon_{C\otimes D}=\varepsilon_C\otimes\varepsilon_D$.
\medskip
{\sl I.3.2. Comodules\/}
\medskip
Let $C$ be any co-unital coalgebra. A {\sl left $C$-comodule\/} is a $k$-vector space $M$ together with a map $\Phi:M\to C\otimes M$ such that the following diagrams commute~:
\diagramme{
\xymatrix{
C\otimes C\otimes M 
	&C\otimes M \ar[l]_{\Delta\otimes I}	&&k\otimes M  	& C\otimes M\ar[l]_{\varepsilon\otimes I}
 \\
C\otimes M \ar[u]_{I\otimes \Phi}	& M \ar[l]_\Phi\ar[u]_\Phi&&&C\ar[u]_\Phi\ar[ul]_\sim}
}
The notion of {\sl right $C$-comodule\/} is defined similarly. A linear subspace $N$ of a left $C$-comodule $M$ is called a {\sl subcomodule\/} if $\Phi(C)\subset C\otimes N$. The intersection of all left subcomodules of $M$ containing a subset $P$ is again a subcomodule, called the {\sl left subcomodule generated by $P$\/}.
\smallskip
It will be convenient to use again Sweedler's notation~:
$$\Phi(m)=\sum_{(m)}m_1\otimes m_0,$$
with $m_0\in M$ and $m_1\in C$. Comodule property reads in Sweedler's notation~:
$$(\Phi\otimes I)\circ\Phi(m)=\sum_{(x)}m_{1:1}\otimes m_{1:2}\otimes m_0=
	\sum_{(m)}m_1\otimes m_{0:1}\otimes m_{0:0}
	=(I\otimes\Delta)\circ\Phi(m),$$
We shall sometimes write the iterated coproduct as~:
$$\sum_{(m)}m_1\otimes m_2\otimes m_0.$$
For a right comodule we have a similar behaviour with $m_0$ on the left. The notion of comodule is dual to the notion of module in the sense that any left (resp. right) $C$-comodule $M$ admits a right (resp. left) $C^*$-module structure. To be precise suppose for the moment that $\Phi$ is any linear map from $M$ to $M\otimes C$, and define $\alpha_{\Phi}:C^*\otimes M\to M$ as the following composition~:
\diagramme{
\xymatrix{C^*\otimes M\ar[r]^{I\otimes\Phi}	&C^*\otimes M\otimes C\
		\ar[r]^{\tau\otimes I}		&M\otimes C^*\otimes C\
		\ar[rr]^{I\otimes <-,->}		&&M\otimes k
		\ar[r]^\sim			&M.\\}
}
Then~:
\prop{I.3.3}
$(M,\Phi)$ is a right (resp. left) $C$-comodule if and only is $(M,\alpha_\Phi)$ is a left (resp. right) $C^*$-module.
\dem
See Sweedler \cite {Sw} section 2.1.
\qed
Note that the duality property is not perfect~: if the linear dual of a coalgebra is always an algebra, the linear dual of an algebra is not in general a coalgebra. However the {\sl restricted dual\/} $A^\circ$ of an algebra $A$ is a coalgebra. It is defined as the space of linear forms on $A$ vanishing on some finite-codimensional ideal. Along the same lines, for a coalgebra $C$ the only left $C^*$-modules related to a right $C$-comodule structure via proposition I.3.3 are the {\sl rational\/} left $C^*$-modules, i.e. those left modules such that the linear map~:
$$\eqalign{\rho : M	&\longrightarrow \Cal L(C^*,M)	\cr
		m	&\longmapsto(x\mapsto x.m)		\cr}$$
has image included in $M\otimes C$ via the embedding~:
$$\eqalign{j:M\otimes C	&\longrightarrow \Cal L(C^*,M)	\cr
	m\otimes c	&\longmapsto (x\mapsto <x,c>m).	\cr}$$
See \cite {Sw} for details. We come now to the fundamental theorem of comodule structure theory~:
\th{I.3.4}
Let $M$ be a left comodule over a coalgebra $C$. For any element $m\in M$ the subcomodule generated by $m$ is finite-dimensional.
\dem
There is a finite collection $(c_i)_{i=1\cdots s}$ of linearly independant elements of $C$ and a collection $(m_i)_{i=1\cdots s}$ of elements of $M$ such that~
$$\Phi(m)=\sum_{i=1}^s m_i\otimes c_i.$$
Let $N$ be the linear subspace of $M$ generated by $\uple ms$. Let us show that $N$ is a left subcomodule of $M$. First note that thanks to the co-unit axiom we have~:
$$m=(I\otimes \varepsilon)\circ \Phi(m)=\sum_{j=1}^s \varepsilon(c_j)m_j,$$
hence $m\in N$. On the other hand, considering linear forms $(f_i)_{i=1\cdots s}$ on $C$ such that $f_i(c_j)=\delta_i^j$ we have~:
$$\eqalign{\Phi(m_i)	&=(I\otimes I\otimes f_i)
	\bigl(\Phi(m_i)\otimes c_i\bigr)	\cr
	&=(I\otimes I\otimes f_i)
		\Bigl(\sum_{j=1}^s\Phi(m_j)\otimes c_j\Bigr) \cr
	&=(I\otimes I\otimes f_i)\circ (\Phi\otimes I)\circ \Phi(m)	\cr
	&=(I\otimes I\otimes f_i)\circ (I\otimes \Delta)\circ \Phi(m)	\cr
	&=\sum_{j=1}^s m_j\otimes (I\otimes f_i)(\Delta c_j),	\cr}$$
whence $\Phi(m_i)\in N\otimes C$, which proves the theorem.
\qed
\cor{I.3.5}
Let $M$ be a left comodule over a coalgebra $C$. Then any left subcomodule of $M$ generated by a finite set is finite-dimensional.
\dem
remark that if $P=\{m_1,\ldots ,m_n\}$, the left subcomodule generated by $P$ is the sum of the left comodules generated by the $m_j$'s, and then apply theorem I.3.4.
\qed
{\sl I.3.3. Structure of coalgebras\/}
\medskip
Let $C$ be a coalgebra. Any intersection of subcoalgebras is a subcoalgebra. To see this consider any family $(D_\alpha)_{\alpha\in \Lambda}$ of subcoalgebras of $C$. Then $I:=\sum_\alpha D_\alpha^\perp$ is a two-sided ideal of $C^*$, as a sum of two-sided ideals. Hence $I^\perp$ is a subcoalgebra according to proposition I.3.2. But $I^\perp$ is indeed the intersection of the subcoalgebras $D_\alpha$.
\ssq
In particular, the intersection of all subcoalgebras containing a given subset $P$ of $C$ will be called the subcoalgebra generated by $P$. We can now state the fundamental theorem of coalgebra theory~:
\th{I.3.6}
Let $C$ be a coalgebra. Then the subcoalgebra generated by one single element $x$ is finite-dimensional.
\dem
The coalgebra $C$ is a left comodule over itself. Let $N$ be the left subcomodule generated by $x$. According to theorem I.3.4, $N$ is finite-dimensional. Then $N^\perp$ has finite codimension, equal to $\mop{dim}N$. It is a left ideal thanks to proposition I.3.1. The quotient space $E=C^*/N^\perp$ is a finite-dimensional left module over $C^\perp$. Let $K$ be the annihilator of this left module. As kernel of the associated representation $\rho:C^*\to\mop{End}E$ it has clearly finite codimension, and it is a two-sided ideal.
\ssq
Now $K^\perp$ is a subcoalgebra according to proposition I.3.2. Moreover it is finite-dimensional, as $\mop{dim}K^\perp=\mop{codim}K^{\perp\perp}\le\mop{codim}K$. Finally $K\subset N^\perp$ implies that $N^{\perp\perp}\subset K^\perp$. A fortiori $N\subset K^\perp$, so $x$ belongs to $K^\perp$. The subcoalgebra generated by $x$ is then included in the finite-dimensional subcoalgebra $K^\perp$, which proves the theorem.
\qed
A coalgebra $C$ is said to be {\sl irreducible\/} if two nonzero subcoalgebras of $C$ have always nonzero intersection. A {\sl simple\/} coalgebra is a coalgebra which does not contain any proper subcoalgebra. A coalgebra $C$ will be called {\sl pointed\/} if any simple subcoalgebra of $C$ is one-dimensional.
\lemme{I.3.7}
Any coalgebra $C$ contains a simple subcoalgebra.
\dem
According to theorem I.3.6 we may suppose that $C$ is finite-dimensional, and the lemma is immediate in this case.
\qed 
\prop{I.3.8}
A coalgebra $C$ is irreducible if and only if it contains a unique simple subcoalgebra.
\dem
Suppose $C$ irreducible, and suppose that $D_1$ and $D_2$ are two simple subcoalgebras. The intersection $D_1\cap D_2$ is nonzero, and hence, by simplicity, $D_1=D_2$. Conversely suppose that $E$ is the only simple subcoalgebra of $C$, and let $D$ any subcoalgebra. According to lemma I.3.7 we have $E\subset D$, hence $E$ is included in any intersection of subcoalgebras, which proves that $C$ is irreducible.
\qed
\lemme{I.3.9}
Let $(C_\alpha)_{\alpha\in\Lambda}$ a family of subcoalgebras of a coalgebra $C$ such that $C$ is the direct sum of the $C_\alpha$'s. Then for any subcoalgebra $D$ we have~:
$$D=\bigoplus_{\alpha\in\Lambda}D\cap C_\alpha.$$
\dem
The sum is indeed direct and included in $D$. To prove the reverse inclusion, pick any $y$ in $D$ and decompose it inside $C$~:
$$y=\sum_{\alpha\in\Lambda}y_\alpha$$
with $y_\alpha\in C_\alpha$ (finite sum). Consider for any $\gamma\in\Lambda$ the linear form $f_\gamma$ defined by~:
$$\eqalign{f_\gamma\restr{C_\alpha}	&=\varepsilon\restr{C_\alpha}
						\hbox{ si }\alpha=\gamma\cr
					&=0\hbox{ si }\alpha\not=\gamma.\cr}$$
Then $f_\gamma(y)=\varepsilon (y_\gamma)$. Now we have~:
$$\eqalign{(I\otimes f_\gamma)\circ\Delta (y)
	&=\sum_\alpha (I\otimes f_\gamma)\circ\Delta (y_\alpha) 	\cr
	&=\sum_\alpha\sum_{(y_\alpha)}(y_\alpha)_1f_\gamma\bigl((y_\alpha)_2\bigr)\cr
	&=\sum_{(y_\gamma)}(y_\gamma)_1\varepsilon\bigl((y_\gamma)_2\bigr)	\cr
	&=(I\otimes\varepsilon)\circ\Delta (y_\gamma)	\cr
	&=y_\gamma.	\cr}$$
This shows that $y_\gamma$ in in $D$, which proves the lemma.
\qed
Let us define the {\sl coradical\/} of a coalgebra $C$ as the sum $R$ of its simple subcoalgebras. As indicated by the terminology this notion is dual to the notion of Jacobson radical of an algebra~: cf. proposition I.3.11 below. 
\prop{I.3.10}
Let $R$ be the coradical of a coalgebra $C$. Then for any subcoalgebra $D$ the coradical $R_D$ of $D$ is equal to $R\cap D$.
\dem
Any simple subcoalgebra of $D$ is a simple subcoalgebra of $C$, so $R_D\subset R\cap D$. Conversely by lemma I.3.9 $R\cap D$ is a direct sum of simple subcoalgebras of $D$, hence $R\cap D\subset R_D$.
\qed
\prop{I.3.11}
If $C$ is a finite-dimensional coalgebra with coradical $R$, then $R^\perp$ is the Jacobson radical of the algebra $C^*$.
\dem
If $S$ is a simple subcoalgebra of $C$ it is clear from dimension considerations that $S^\perp$ is a maximal two-sided ideal of $C^*$. Conversely any maximal two-sided ideal of $C^*$ is the orthogonal of a simple subcoalgebra, so $R^\perp$ is indeed the intersection of all maximal two-sided ideals of $C^*$. Finally (lemma I.2.8), maximal two-sided and primitive ideals of a finite-dimensional algebra coincide.
\qed
\medskip
{\sl I.3.4. The wedge\/}
\medskip
Let $C$ be a coalgebra and $X,Y$ two linear subspaces of $C$. We define~:
$$X\wedge Y=\{x\in C,\ \Delta x\in X\otimes C+C\otimes Y\}.$$
We define as well inductively~:
$$\wedge^0 X=\{0\},\hskip 12mm \wedge^nX=(\wedge^{n-1}X)\wedge X.$$
The alternative definition in terms of the algebra structure on $C^*$ is often more manageable, and its verification is straightforward~:
$$X\wedge Y=(X^\perp Y^\perp)^\perp.$$
Here is a first application of this definition~:
\prop{I.3.12}
1) The wedge is associative : $(X\wedge Y)\wedge Z=X\wedge(Y\wedge Z)$.

2) If $X$ is a left coideal then $\{0\}\wedge X=X$.

3) If $X$ is a left coideal and $Y$ is a right coideal, then $X\wedge Y$ is a subcoalgebra of $C$.

4) The wedge of two subcoalgebras is a subcoalgebra.

5) If $X\subset X'$ and $Y\subset Y'$, then $X\wedge Y\subset X'\wedge Y'$.
\dem
According to the definition we have~:
$$(X\wedge Y)\wedge Z=X\wedge(Y\wedge Z)=(X^\perp Y^\perp Z^\perp)^\perp.$$
Let $X$ (resp. $Y$) be a left (resp. right) coideal. According to proposition I.3.1, $X^\perp$ is a left ideal of $C^*$ and $Y^\perp$ is a right ideal. The product $(X^\perp Y^\perp)$ is then a two-sided ideal. Second and third assertions follow by applying proposition I.3.1 again, and by noticing that~:
$$\{0\}\wedge X=(C^*.X^\perp)^\perp=X^{\perp\perp}=X.$$
4) is an immediate consequence of 3), and 5) is clear.
\qed
The wedge admits the following comodule version~: let $M$ be a right $C$-comodule with coaction $\Phi:M\to M\otimes C$. Let $N$ be a subspace of $M$ and $X$ be a subspace of $C$. We define~:
$$N\wedge X=\{x\in M,\ \Phi x\in N\otimes C+C\otimes X\}.$$
One can check that if $X$ is a right coideal the $N\wedge X$ is a subcomodule, and if $N$ is a subcomodule then $N\subset N\wedge X$.
\prop{I.3.13} 
Let $R$ be the coradical of a coalgebra $C$, and let $M$ be a right $C$-comodule. Then for $\{0\}\subset M$ we have~:
$$\bigcup_n\{0\}\wedge(\wedge^nM)=M.$$
\dem
Let $x\in M$, and let $N$ be a finite-dimensional subcomodule containing $x$ (which exists thanks to theorem I.3.4). If $\Phi$ denotes the coaction, then $\Phi(N)\subset N\otimes X$ where $X$ is a finite-dimensional subspace of $C$. Let $D$ be a finite-dimensional subcoalgebra containing $X$ (which exists thanks to theorem I.3.6). It is clear that $N$ is a right $D$-comodule.
\ssq
Applying Proposition I.3.10, the coradical of $D$ is $R_0=R\cap D$. Proposition I.3.11 says that $R_0^\perp$ is the Jacobson radical of $D^*$, which is nilpotent by corollary I.2.7~: there exists a positive integer $n$ such that $(R_0^\perp)^n=\{0\}$. Dualizing we get that $\wedge^n_DR_0=D$, where the subscribe $D$ reminds with respect to which coalgebra the wedge operation is performed. Clearly we have~:
$$\wedge^n_DR_0\subset \wedge^n_CR_0\subset \wedge^n_CR,$$
the second inclusion coming from assertion 5) of proposition I.3.12. We have then~:
$$D\subset\wedge^nR$$
(we have dropped the subscribe ``$C$'' here), hence the inclusions~:
$$N\subset \{0\}\wedge D\subset \{0\}\wedge(\wedge^nR).$$
The initial element $x$ belongs then to $\{0\}\wedge(\wedge^nR)$ for some $n$, which proves the assertion.
\qed
\medskip
{\sl I.3.5. The coradical filtration\/}
\medskip
Let $C$ be a coalgebra with coradical $R$. We consider for any integer $i\ge 0$~:
$$C^i=\wedge^{i+1}R.$$
The following proposition is an immediate consequence of propositions I.3.12 and I.3.13~:
\prop{I.3.14}
$(C^i)_{i\ge 0}$ is an increasing sequence of subcoalgebras of $C$, and we have~:
$$C=\bigcup_{i\ge 0}C^i.$$
\ndem
The coalgebra $C$ is then endowed with an increasing filtration by subcoalgebras~: its {\sl coradical filtration\/}.
\prop{I.3.15}
The coproduct is compatible with the coradical filtration, in the sense that the following inclusion holds~:
$$\Delta(C^n)\subset\sum_{i=0}^n C^i\otimes C^{n-i}.$$
\dem
For any $i\in\{0,\ldots ,n+1\}$ we have~:
$$\wedge^{n+1}R=(\wedge^iR)\wedge(\wedge^{n-i+1}R).$$
This is immediate for $i=1,\ldots ,n$, and comes from proposition I.3.12 assertion 2) for $i=0$ and $i=n+1$. We have then (setting $C^{-1}=\{0\}$)~:
$$\Delta C^n\subset\bigcap_{i=0}^{n+1}(C\otimes C^{n-i}+C^{i-1}\otimes C).$$
The right-hand side $(RHS)$ is contained in $C^n\otimes C^n$. Choose any supplementary subspace $D_i$ of $C^{i-1}$ inside $C^i$~: so $C^0=D_0$, and $C^i=D_0\oplus\cdots\oplus D_i$. We have~:
$$\eqalign{(RHS)	&=\bigcap_{i=0}^{n+1}\bigoplus_
{r\le i-1\hbox{ or } s\le n-i}D_r\otimes D_s	\cr
			&=\bigoplus_{r+s\le n}D_r\otimes D_s	\cr
			&=\sum_{i=0}^n C^i\otimes C^{n-i},	\cr}$$
which proves the result.
\qed

\alinea{I.4. Convolution product}
Let $A$ be an algebra and $C$ be a coalgebra (over the same field $k$). Then there is an associative product on $\Cal L(C,A)$ called the {\sl convolution product\/}. It is given by~:
$$\varphi*\psi=m_{A}\circ(\varphi\otimes\psi)\circ\Delta_{C}.$$
In Sweedler's notation it reads~:
$$\varphi*\psi(x)=\sum_{(x)}\varphi(x_1)\psi(x_2).$$
The associativity is a direct consequence of both associativity of $A$ and coassociativity of $C$. We shall study this product more thoroughly in the next chapters.
\alinea{I.5. Bialgebras and Hopf algebras}
A (unital and co-unital) {\sl bialgebra\/} is a vector space $\Cal H$ endowed with a structure of unital algebra $(m,u)$ and a structure of co-unital coalgebra $(\Delta,\varepsilon)$ which are compatible. The compatibility requirement is that $\Delta$ is an algebra morphism (or equivalently that $m$ is a coalgebra morphism), $\varepsilon$ is an algebra morphism and $u$ is a coalgebra morphism. It is expressed by the commutativity of the three following diagrams~:
\diagramme{
\xymatrix{\Cal H\otimes\Cal H\otimes\Cal H\otimes\Cal H
\ar[rr]^{\tau_{23}}	&&\Cal H\otimes\Cal H\otimes\Cal H\otimes\Cal H
\ar[d]^{m\otimes m}	\\
\Cal H\otimes \Cal H	\ar[u]_{\Delta\otimes\Delta} \ar[r]_m
&\Cal H \ar[r]_{\Delta}	&\Cal H\otimes\Cal H}
}
\diagramme{
\xymatrix{\Cal H\otimes\Cal H \ar[d]^m \ar[r]^{\varepsilon\otimes\varepsilon}
	&k\otimes k	\ar[d]^\sim	&&&\Cal H\otimes \Cal H	&k\otimes k \ar[l]_{u\otimes u}\\
\Cal H\ar[r]^\varepsilon &k	&&&\Cal H \ar[u]_\Delta	&k\ar[l]_u \ar[u]_\sim}
}
A {\sl Hopf algebra\/} is a bialgebra $\Cal H$ together with a linear map $S:\Cal H\to \Cal H$ called the {\sl antipode\/}, such that the following diagram commutes~:
\diagramme{
\xymatrix{&\Cal H\otimes\Cal H	\ar[rr]^{S\otimes I}
				&&\Cal H\otimes \Cal H\ar[dr]^{m}	& \\
\Cal H\ar[rr]^\varepsilon \ar[dr]^\Delta \ar[ur]^\Delta
				&&	k\ar[rr]^u	&&\Cal H\\
&\Cal H\otimes\Cal H	\ar[rr]^{I\otimes S}
				&&\Cal H\otimes \Cal H\ar[ur]^{m}	&}
}
In Sweedler's notation it reads~:
$$\sum_{(x)}S(x_1)x_2=\sum_{(x)}x_1S(x_2)=(u\circ\varepsilon)(x).$$
In other words the antipode is an inverse of the identity $I$ for the convolution product on $\Cal L(H,H)$. The unit for the convolution is the map $u\circ\varepsilon$.
\ssq
Let $\Cal H$ be a bialgebra. A {\sl primitive element\/} in $\Cal H$ is an element $x$ such that $\Delta x=x\otimes 1+1\otimes x$. A {\sl grouplike element\/} is a nonzero element $x$ such that $\Delta x=x\otimes x$. Note that grouplike elements make sense in any coalgebra.
\alinea{I.6. Examples}
{\sl I.6.1. The Hopf algebra of a group}\msq
Let $G$ be a group, and let $kG$ be the group algebra (over the field $k$). It is by definition the vector space freely generated by the elements of $G$~: the product of $G$ extends uniquely to a bilinear map from $kG\times kG$ into $kG$, hence a multiplication $m:kG\otimes kG\to kG$, which is associative. The neutral element of $G$ gives the unit for $m$.
\smallskip
The space $kG$ is also endowed with a co-unital coalgebra structure, given by~:
$$\Delta (\sum \lambda_ig_i)=\sum \lambda_i.g_i\otimes g_i$$
and~:
$$\varepsilon (\sum \lambda_ig_i)=\sum \lambda_i.$$
This defines the {\sl coalgebra of the set $G$\/}~: it does not take into account the extra group structure on $G$, as the algebra structure does.
\prop{I.6.1}  
The vector space $kG$ endowed with the algebra and coalgebra structures defined above is a Hopf algebra. The antipode is given by~:
$$S(g)=g\inver, g\in G.$$
\dem
The compatibility of the product and the coproduct is an immediate consequence of the following computation~: for any $g,h\in G$ we have~:
$$\Delta(gh)=gh\otimes gh=(g\otimes g)(h\otimes h)=\Delta g.\Delta h.$$
Now $m(S\otimes I)\Delta(g)=g\inver g=e$ and similarly for $m(I\otimes S)\Delta(g)$. But $e=u\circ\varepsilon(g)$ for any $g\in G$, so map $S$ is indeed the antipode.
\qed
{\sl Remark\/}~: if $G$ were only a semigroup, the same construction would lead to a bialgebra structure on $kG$~: the Hopf algebra structure (i.e. the existence of an antipode) reflects the group structure (the existence of the inverse). We have $S^2=I$ in this case, but involutivity of the antipode is not true for general Hopf algebras.
\medskip
{\sl I.6.2. Tensor algebras\/}\msq
There is a natural structure of cocommutative Hopf algebra on the tensor algebra $T(V)$ of any vector space $V$. Namely we define the coproduct $\Delta$ as the unique algebra morphism from $T(V)$ into $T(V)\otimes T(V)$ such that~:
$$\Delta(1)=1\otimes 1,\hskip 12mm \Delta(x)=x\otimes 1+1\otimes x, \ x\in V.$$
We define the co-unit as the algebra morphism such that $\varepsilon(1)=1$ and $\varepsilon\restr V=0$ This endows $T(V)$ with a cocommutative bialgebra structure. We claim that the principal anti-automorphism~:
$$S(x_1\otimes\cdots\otimes x_n)=(-1)^nx_n\otimes\cdots\otimes x_1$$
verifies the axioms of an antipode, so that $T(V)$ is indeed a Hopf algebra. For $x\in V$ we have $S(x)=-x$, hence $S*I(x)=I*S(x)=0$. As $V$ generates $T(V)$ as an algebra it is easy to conclude.
\medskip
{\sl I.6.3. Enveloping algebras\/}\msq
Let $\g g$ a Lie algebra. The universal enveloping algebra is the quotient of
the tensor algebra $T(\g g)$ by the ideal $J$ generated by $x\otimes
y-y\otimes x-[x,y], \ x,y\in\g g$.
\lemme{I.6.2}
$J$ is a Hopf ideal, i.e. $\Delta(J)\subset\Cal H\otimes J+J\otimes\Cal H$ and $S(\Cal H)\subset\Cal H$.
\dem
The last assertion is immediate. The first comes easily from the fact that the ideal $J$ is generated by primitive elements (according to proposition I.7.3 below)~: indeed any ideal generated by primitive elements is a Hopf ideal (very easy and left to the reader).
\qed
The quotient of a Hopf algebra by a Hopf ideal is a Hopf algebra. Hence the universal enveloping algebra $\Cal U(\g g)$ is a cocommutative Hopf algebra.  
\alinea{I.7. Some properties of Hopf algebras}
We summarize in the proposition below the main properties of the antipode in a Hopf algebra~:
\prop{I.7.1 \rm(cf. \cite{Sw} proposition 4.0.1)}
Let $\Cal H$ be a Hopf algebra with multiplication $m$, comultiplication $\Delta$, unit $u:1\mapsto \un$, co-unit $\varepsilon$ and antipode $S$. Then~:
\smallskip
1) $S\circ u=u$ and $\varepsilon\circ S=\varepsilon$.
\smallskip
2) $S$ is an algebra antimorphism and a coalgebra antimorphism, i.e. if $\tau$ denotes the flip we have~:
$$m\circ (S\otimes S)\circ\tau=S\circ m,\hskip 12mm \tau\circ(S\otimes S)\circ\Delta=\Delta\circ S.$$
3) If $\Cal H$ is commutative or cocommutative, then $S^2=I$.
\dem
We follow closely the proof given by Chr. Kassel here \cite K. Starting from $\varepsilon(\un)=1$ and $\Delta(\un)=\un\otimes\un$, we get~: $\un=m(S\otimes I)\Delta(\un)=S(\un)\un$, so $S(\un)=\un$. We have for any $x\in\Cal H$~:
$$\varepsilon\circ u\circ\varepsilon (x)=\varepsilon (\varepsilon(x).\un)=\varepsilon(x).$$
Then,
$$\eqalign{\varepsilon (x)	
	&=\varepsilon\bigl(m(S\otimes I)\Delta\bigr)(x)	\cr
	&=\varepsilon\bigl(\sum_{(x)}S(x_1)x_2\bigl)	\cr
	&=\sum_{(x)}(\varepsilon\circ S)(x_1)\varepsilon(x_2).}$$
On the other hand,
$$S(x)=S*(u\circ\varepsilon)(x)=\sum_{(x)}S(x_1)(u\circ\varepsilon)(x_2).$$
Hence,
$$(\varepsilon\circ S)(x)=\sum_{(x)}(\varepsilon\circ S)(x_1)\varepsilon (x_2).$$
This proves assertion 1). Now consider $m,N,P\in\cal L(\Cal H\otimes\Cal H,\Cal H)$ defined by~:
$$m(x\otimes y)=xy,\hskip 12mm N(x\otimes y)=S(y)S(x), \hskip 12mm
P(x\otimes y)=S(xy).$$
Considering the convolution product $\wt *$ on $\cal L(\Cal H\otimes\Cal H,\Cal H)$ we shall prove~:
\lemme{I.7.2}
$$P\wt *m=u\circ\varepsilon_{\cal H\otimes\cal H}=m\wt *N.$$
\dem
We compute, with Sweedler's notation~:
$$\eqalign{P\wt *m(x\otimes y)	&=\sum_{(x\otimes y)}P\bigl((x\otimes y)_1\bigr)
				m\bigl((x\otimes y)_2\bigr)	\cr
			&=\sum_{(x),(y)}P(x_1\otimes y_1)
				m(x_2\otimes y_2)		\cr
			&=\sum_{(x),(y)}S(x_1y_1)x_2y_2		\cr
			&=(S*I)(xy)				\cr
			&=(u\circ\varepsilon) (xy)		\cr
			&=(u\circ\varepsilon_{\Cal H\otimes\Cal H})
					(x\otimes y),		\cr}$$
and~:
$$\eqalign{m\wt *N(x\otimes y)	&=\sum_{(x),(y)}m(x_1\otimes y_1)
					N(x_2\otimes y_2)	\cr
			&=\sum_{(x),(y)}x_1y_1S(y_2)S(x_2)	\cr
			&=\sum_{(x)}x_1(u\circ\varepsilon)(y)S(x_2)\cr
			&=(u\circ\varepsilon)(x)(u\circ\varepsilon)(y)\cr
			&=(u\circ\varepsilon_{\Cal H\otimes\Cal H})
					(x\otimes y),\cr}$$
which proves the lemma.
\qed
{\sl End of proof of Proposition I.7.1\/}~: as $u\circ\varepsilon_{\Cal H\otimes\Cal H}$ is the unit element for $\wt *$ lemma proves that both $P$ and $N$ are inverse of $m$ for this convolution. Hence~:
$$P=P\wt *(m\wt * N)=(P\wt * m)\wt * N=N,$$
which proves the first part of assertion 2). The second part is proved similarly using convolution in $\Cal L(\Cal H,\Cal H\otimes H)$ with $m$ replaced with $\Delta$, $N$ replaced with $\tau\circ(S\otimes S)\circ\Delta$ and $P$ replaced with $\Delta\circ S$. Indeed,
$$\eqalign{\bigl((\Delta\circ S)\wt * \Delta\bigr)(x)
	&=(m\otimes m)\circ\tau_{23}\circ
		\bigl((\Delta\circ S)\otimes \Delta\bigr)\circ\Delta(x)	\cr
	&=(m\otimes m)\circ\tau_{23}\circ(\Delta\otimes\Delta)\circ(S\otimes I)
		\circ\Delta(x)	\cr
	&=\Delta\circ m\circ (S\otimes I)\circ\Delta(x)	\cr
	&=\Delta\bigl(u\circ\varepsilon(x)\bigr)	\cr
	&=u\circ\varepsilon(x)\otimes u\circ\varepsilon(x),\cr}$$
and~:
$$\eqalign{\hskip -12mm\Delta\wt *\bigl(\tau\circ(S\otimes S)\circ\Delta\bigr)(x)
	&=(m\otimes m)\circ\tau_{23}\circ
		\Bigl(\Delta\otimes\bigl(\tau\circ(S\otimes S)\circ\Delta\bigr)
		\Bigr)\circ \Delta(x)	\cr
	&=(m\otimes m)\circ\tau_{23}\circ\tau_{34}\circ
		(I\otimes I\otimes S\otimes S)\circ (\Delta\otimes\Delta)
		\circ\Delta(x)	\cr
	&=(m\otimes m)\circ\tau_{23}\circ\tau_{34}\circ
		(I\otimes I\otimes S\otimes S)\circ (I\otimes\Delta\otimes I)
		\circ(\Delta\otimes I)\circ\Delta(x)	\cr
	&= (m\otimes m)\circ\tau_{23}\circ\tau_{34}\circ 
		(I\otimes I\otimes S\otimes S)\circ(I\otimes\Delta\otimes I)
		\Bigl(\sum_{(x)}x_1\otimes x_2\otimes x_3\Bigr)	\cr
	&=(m\otimes m)\circ\tau_{23}\circ\tau_{34}\circ 
		\Bigl(\sum_{(x)}x_1\otimes x_{2:1}\otimes Sx_{2:2}
		\otimes Sx_3\Bigr)	\cr
	&=\sum_{(x)}x_1Sx_3\otimes x_{2:1}Sx_{2:2}	\cr
	&=\sum_{(x)}x_1Sx_3\otimes(u\circ\varepsilon)(x_2)	\cr
	&=\Bigl(\sum_{(x)}x_1Sx_2\Bigr)\otimes(u\circ\varepsilon)(x)	\cr
	&=u\circ\varepsilon(x)\otimes u\circ\varepsilon(x).\cr}$$

If $\Cal H$ is commutative then~:
$$u\circ\varepsilon(x)=\sum_{(x)}x_2S(x_1).$$
If $\Cal H$ is cocommutative, this is again true, as~:
$$u\circ\varepsilon(x)=(I*S)(x)=\sum_{(x)}x_1S(x_2).$$
Suppose that $\Cal H$ is commutative or cocommutative. Then using assertion 2) we can compute~:
$$\eqalign{(S\circ S)*S(x)	&=\sum_{(x)}(S\circ S)(x_1)S(x_2)	\cr
			&=S\bigl(\sum_{(x)}x_2S(x_1)\bigr)	\cr
			&=S\bigl(u\circ\varepsilon(x)\bigr)	\cr
			&=(u\circ\varepsilon)(x).		\cr}$$
Thus $S\circ S$ is a left inverse for $S$ (for the convolution), and then $S\circ S=I$. This ends the proof of Proposition I.7.1.
\qed
\prop{I.7.3}
1). If $x$ is a primitive element then $S(x)=-x$.
\smallskip
2). The linear subspace $\mop{Prim}\Cal H$ of primitive elements in $\Cal H$ is a Lie algebra.
\dem
If $x$ is primitive, then $(\varepsilon\otimes\varepsilon)\circ\Delta (x)=2\varepsilon(x)$. On the other hand, $(\varepsilon\otimes\varepsilon)\circ\Delta (x)=\varepsilon(x)$, so $\varepsilon(x)=0$. Then~:
$$0=(u\circ\varepsilon)(x)=m(S\otimes I)\Delta(x)=S(x)-x. $$
Now let $x$ and $y$ be primitive elements of $\Cal H$. Then we can easily compute~:
$$\eqalign{\Delta(xy-yx)	&=(x\otimes \un+\un\otimes x)(y\otimes \un
+\un\otimes y)-(y\otimes \un+\un\otimes y)(x\otimes \un+\un\otimes x)	\cr
				&=(xy-yx)\otimes \un+\un\otimes (xy+yx)
+x\otimes y+y\otimes x-y\otimes x-x\otimes y	\cr
				&=(xy-yx)\otimes \un+\un\otimes (xy-yx).\cr}$$
\qed
\paragraphe{II. Convolution product and regularization}
\alinea{II.1. Connected graded bialgebras}
Let $k$ be a field with vanishing characteristic. We shall denote by $k[[t]]$ the ring of formal series on $k$, and by $k[t\inver,t]]$ the field of Laurent series on $k$. A {\sl graded Hopf algebra\/} on $k$ is a graded $k$-vector space~:
$$\Cal H=\bigoplus_{n\ge 0}\Cal H_n$$
endowed with a product $m:\Cal H\otimes \Cal H\to\Cal H$, a coproduct $\Delta:\Cal H\to\Cal H\otimes \Cal H$, a unity $u:k\to\Cal H$, a co-unity $\varepsilon:\Cal H\to k$ and an antipode $S:\Cal H\to\Cal H$ fulfilling the usual axioms of a Hopf algebra, and such that~:
$$\eqalign{\Cal H_p.\Cal H_q	&\subset \Cal H_{p+q}	\cr
		\Delta(\Cal H_n)	&\subset \bigoplus_{p+q=n}\Cal H_p\otimes\Cal H_q. \cr
	S(\Cal H_n)	&\subset\Cal H_n	\cr}$$
If we do not ask for the existence of an antipode $\Cal H$ we get the definition of a {\sl graded bialgebra\/}. In a graded bialgebra $\Cal H$ we shall consider the increasing filtration~:
$$\Cal H^n=\bigoplus_{p=0}^n\Cal H_p.$$
Suppose moreover that $\Cal H$ is {\sl connected\/}, i.e. $\Cal H_{0}$ is one-dimensional. Then we have~:
$$\mop{Ker}\varepsilon=\bigoplus_{n\ge 1}\Cal H_n.$$
\prop{II.1.1}
For any $x\in\Cal H^n, n\ge 1$ we can write~:
$$\Delta x=x\otimes\un+\un\otimes x+\wt\Delta x,\hskip 12mm 
\wt\Delta x\in\bigoplus_{p+q=n,p\not= 0,q\not= 0}\Cal H_p\otimes\Cal H_q.$$
The map $\wt\Delta$ is coassociative on $\mop{Ker}\varepsilon$ and $\wt\Delta_k=(I^{\otimes k-1}\otimes\wt\Delta)(I^{\otimes k-2}\otimes\wt\Delta)...\wt\Delta$ sends $\Cal H^n$ into $(\Cal H^{n-k})^{\otimes k+1}$.
\dem
Thanks to connectedness we clearly can write~:
$$\Delta x=a(x\otimes 1)+b(1\otimes x)+\wt\Delta x$$
with $a,b\in k$ and $\wt\Delta x\in\mop{Ker}\varepsilon\otimes\mop{Ker}\varepsilon$. The co-unity property then tells us that, with $k\otimes \Cal H$ and $\Cal H\otimes k$ canonically identified with $\Cal H$~:
$$\eqalign{x	&=(\varepsilon\otimes I)(\Delta x)=bx	\cr
	x	&=(I\otimes\varepsilon)(\Delta x)=ax,	\cr}$$
hence $a=b=1$. We shall use the following two variants of Sweedler's notation~:
$$\eqalign{\Delta x	&=\sum_{(x)}x_1\otimes x_2,	\cr
	\wt\Delta x	&=\sum_{(x)}x'\otimes x'',	\cr}$$
the second being relevant only for $x\in\mop{Ker}\varepsilon$. if $x$ is homogeneous of degree $n$ we can suppose that the components $x_1,x_2,x',x''$ in the expressions above are homogeneous as well, and we have then $|x_1|+|x_2|=n$ and $|x'|+|x''|=n$ We easily compute~:
$$\eqalign{(\Delta\otimes I)\Delta(x)	&=x\otimes 1\otimes 1 +1\otimes x\otimes 1
						+1\otimes 1\otimes x	\cr
				&\ +\sum_{(x)}
x'\otimes x''\otimes 1 + x'\otimes 1\otimes x'' + 1\otimes x'\otimes x''\cr
				&\ +(\wt\Delta\otimes I)\wt \Delta (x)\cr}$$
and
$$\eqalign{(I\otimes\Delta)\Delta(x)	&=x\otimes 1\otimes 1 +1\otimes x\otimes 1
						+1\otimes 1\otimes x	\cr
				&\ +\sum_{(x)}
x'\otimes x''\otimes 1 + x'\otimes 1\otimes x'' + 1\otimes x'\otimes x''\cr
				&\ +(I\otimes\wt\Delta)\wt \Delta (x),\cr}$$
hence the co-associativity of $\wt\Delta$ comes from the one of $\Delta$. Finally it is easily seen by induction on $k$ that for any $x\in\Cal H^n$ we can write~:
$$\wt\Delta_k(x)=\sum_{x}x^{(1)}\otimes\cdots\otimes x^{(k+1)},$$
with $|x^{(j)}|\ge 1$. The grading imposes~:
$$\sum_{j=1}^{k+1}|x^{(j)}|=n,$$
so the maximum possible for any degree $|x^{(j)}|$ is $n-k$.
\qed
\alinea{II.2. Connected filtered bialgebras}
A {\sl filtered Hopf algebra\/} on $k$ is a $k$-vector space together with an increasing $\Z_+$-indexed filtration~:
$$\Cal H^0\subset\Cal H^1\subset\cdots\subset \Cal H^n\subset\cdots, \bigcup_n \Cal H^n=\Cal H$$
endowed with a product $m:\Cal H\otimes \Cal H\to\Cal H$, a coproduct $\Delta:\Cal H\to\Cal H\otimes \Cal H$, a unit characteristic $u:k\to\Cal H$, a co-unit $\varepsilon:\Cal H\to k$ and an antipode $S:\Cal H\to\Cal H$ fulfilling the usual axioms of a Hopf algebra, and such that~:
$$\eqalign{\Cal H^p.\Cal H^q	&\subset \Cal H^{p+q}	\cr
		\Delta(\Cal H^n)	&\subset \sum_{p+q=n}\Cal H^p\otimes\Cal H^q \cr
	S(\Cal H^n)	&\subset\Cal H^n.	\cr}$$
If we do not ask for the existence of an antipode $\Cal H$ we get the definition of a {\sl filtered bialgebra\/}. For any $x\in\Cal H$ we set~:
$$|x|=\mop{min}\{n\in\N,\ x\in\Cal H^n\}.$$
Any graded bialgebra or Hopf algebra is obviously filtered by the canonical filtration associated to the grading~:
$$\Cal H^n=\bigoplus_{i=0}^n \Cal H_i,$$
and in that case, if $x$ is an homogeneous element, $x$ is of degree $n$ if and only if $|x|=n$. We say that the filtered bialgebra $\Cal H$ is connected if $\Cal H^0$ is one-dimensional. There is an analogue of proposition II.1.1 in the connected filtered case~:
\prop{II.2.1}
For any $x\in\Cal H^n, n\ge 1$ we can write~:
$$\Delta x=x\otimes\un+\un\otimes x+\wt\Delta x,\hskip 12mm 
\wt\Delta x\in\sum_{p+q=n,p\not= 0,q\not= 0}\Cal H^p\otimes\Cal H^q.$$
The map $\wt\Delta$ is coassociative on $\mop{Ker}\varepsilon$ and $\wt\Delta_k=(I^{\otimes k-1}\otimes\wt\Delta)(I^{\otimes k-2}\otimes\wt\Delta)...\wt\Delta$ sends $\Cal H^n$ into $(\Cal H^{n-k})^{\otimes k+1}$.
\dem
Straightforward adaptation of proof of proposition II.1.1.
\qed
The following theorem is due to S. Montgomery \cite {Mo lemma 1.1}.
\th{II.2.2}
Let $\Cal H$ be any pointed Hopf algebra. Then the coradical filtration endows $\Cal H$ with a structure of filtered Hopf algebra. 
\dem
 It only remains to show that for any $n\in\N$ the inclusion $S(H^n)\subset H^n$ holds, and that for any $p,q\in\N$~:
$$\Cal H^p\Cal H^q\subset \Cal H^{p+q},$$
which, together with proposition I.3.15, will prove the result. Recall that a pointed coalgebra is a coalgebra in which any simple subcoalgebra is one-dimensional. In this case any simple subcoalgebra is linearly generated by a unique grouplike element. Any grouplike element $g$ in a Hopf algebra admits an inverse $Sg$, where $S$ is the antipode. It follows that the coradical $\Cal H^0$ of a pointed Hopf algebra $\Cal H$ is a Hopf subalgebra of $\Cal H$, precisely the Hopf algebra of the group  of the grouplike elements of $\Cal H$ (cf. example I.6.1).
\ssq
The proof proceeds by induction~: inclusion $S\Cal H_0\subset \Cal H_0$ obviously holds. Suppose that $S\Cal H^k\subset \Cal H^k$ for all $k\le n-1$. Using the definition of $\Cal H^n$~:
$$\Cal H^n=\Cal H^0\wedge\Cal H^{n-1}=\Cal H^{n-1}\wedge H^0$$
and the formula~:
$$Sx=\sum_{(x)}Sx_2\otimes Sx_1,$$
(cf. proposition I.7.1) we deduce the inclusion $S\Cal H^n\subset \Cal H^n$. Now, suppose that the inclusion $\Cal H^k\Cal H^0\subset\Cal H^k$ holds for $k\le n-1$ (its is obviously the case for $k=0$). Then we have~:
$$\eqalign{\Delta (\Cal H^n\Cal H^0)	&\subset (\Cal H^0\otimes \Cal H+\Cal H\otimes \Cal H^{n-1})(\Cal H^0\otimes\Cal H^0)	\cr
	&\subset \Cal H^0\otimes \Cal H+\Cal H\otimes \Cal H^{n-1}\Cal H^0\cr
	&\subset \Cal H^0\otimes \Cal H+\Cal H\otimes \Cal H^{n-1}.\cr}$$
So $\Cal H^n\Cal H^0\subset \Delta\inver(\Cal H^0\otimes \Cal H+\Cal H\otimes \Cal H^{n-1})=\Cal H^n$. Similarly on the other side we have $\Cal H^0\Cal H^n\subset\Cal H^n$ for any $n$. Suppose now that the inclusion~:
$$\Cal H^p\Cal H^q\subset \Cal H^{p+q}$$
holds for any $p,q$ such that $p+q\le n-1$. Choose $p,q$ with $p+q=n$ and compute~:
$$\hskip -4pt\eqalign{\Delta(\Cal H^p\Cal H^q)	&\subset (\Cal H^0\otimes \Cal H^p+\Cal H^p\otimes \Cal H^{p-1})(\Cal H^0\otimes \Cal H^q+\Cal H^q\otimes \Cal H^{q-1})\cr
	&\subset\Cal H^0\Cal H^0\otimes \Cal H^p\Cal H^q
		+\Cal H^p\Cal H^0\otimes \Cal H^{p-1}\Cal H^q
		+\Cal H^0\Cal H^q\otimes \Cal H^p\Cal H^{q-1}
		+\Cal H^p\Cal H^q\otimes \Cal H^{p-1}\Cal H^{q-1}	\cr
	&\subset \Cal H^0\otimes\Cal H+\Cal H\otimes\Cal H^{p+q-1}\cr}$$
thanks to the induction hypothesis and the property already proved when one of the indices is equal to zero. Thus $\Cal H^p\Cal H^q\subset\Cal H^{p+q}$, which finishes the proof of the theorem.
\qed 
{\sl Remark 1}~: The proof used only the property that the coradical is a Hopf subalgebra of $\cal H$. The pointedness of $\Cal H$ implies this property but is not strictly necessary.
\smallskip
{\sl Remark 2}~: the image of $k$ under the unit map $u$ is a one-dimensional simple subcoagebra of $\Cal H$. If $\Cal H$ is an irreducible coalgebra, by proposition I.3.8 it is the unique one, and then the coradical is $\Cal H^0=k.\un$. Any irreducible Hopf algebra is then pointed, and connected with respect to the coradical filtration. 
\alinea{II.3. The convolution product}
\qquad
An important result is that any connected filtered bialgebra is indeed a filtered Hopf algebra, in the sense that the antipode comes for free. We give a proof of this fact as well as a recursive formula for the antipode with the help of the {\sl convolution product\/}~:
\ssq 
Let $\Cal H$ be a (connected filtered) bialgebra, and let $\Cal A$ be any $k$-algebra~: the convolution product on $\Cal L(\Cal H,\Cal A)$ is given by~:
$$\eqalign{\varphi*\psi (x)	&=m_{\Cal A}(\varphi\otimes\psi)\Delta(x)\cr
				&=\sum_{(x)}\varphi(x_1)\psi(x_2).\cr}$$
\prop{II.3.1}
The map $e=u_{\Cal A}\circ\varepsilon$, given by $e(\un)=\un_{\Cal A}$ and $e(x)=0$ for any $x\in\mop{Ker}\varepsilon$, is a unit for the convolution product. Moreover the set $G=\{\varphi\in\Cal L(\Cal H,\Cal A),\ \varphi(\un)=\un_{\Cal A}\}$ endowed with the convolution product is a group.
\dem
The first statement is straightforward. To prove the second let us consider the formal series~:
$$\eqalign{\varphi^{*-1}(x)	&=\bigl(e-(e-\varphi)\bigr)^{*-1}(x)	\cr
				&=\sum_{k\ge 0}(e-\varphi)^{*k}(x).	\cr}$$
Using $(e-\varphi)(\un)=0$ we have immediately $(e-\varphi)^{*k}(\un)=0$, and for any $x\in\mop{Ker}\varepsilon$~:
$$(e-\varphi)^{*k}(x)=m_{\Cal A,k-1}(\varphi\otimes\cdots\otimes\varphi)\wt\Delta_{k-1}(x).$$
When $x\in\Cal H^n$ this expression vanishes then for $k\ge n+1$. The formal series ends up then with a finite number of terms for any $x$, which proves the result.
\qed
\cor{II.3.2}
Any connected filtered bialgebra $\Cal H$ is a filtered Hopf algebra. The antipode is defined by~:
$$S(x)=\sum_{k\ge 0}(u\varepsilon-I)^{*k}(x).$$
It is given by $S(\un)=\un$ and recursively by any of the two formulas for $x\in\mop{Ker}\varepsilon$~:
$$\eqalign{S(x)	&=-x-\sum_{(x)}S(x')x''	\cr
		S(x)	&=-x-\sum_{(x)}x'S(x'').	\cr}$$
\dem
The antipode, when it exists, is the inverse of the identity for the convolution product on $\Cal L(\Cal H,\Cal H)$. One just needs then to apply proposition II.3.1 with $\Cal A=\Cal H$. The two recursive formulas come directly from the two equalities~:
$$m(S\otimes I)\Delta (x)=m(I\otimes S)\Delta (x)=0$$
fulfilled by any $x\in\mop{Ker}\varepsilon$.
\qed
Let $\g g$ be the subspace of $\Cal L(\Cal H,\Cal A)$ formed by the elements $\alpha$ such that $\alpha(\un)=0$. It is clearly a subalgebra of $\Cal L(\Cal H,\Cal A)$ for the convolution product. We have~:
$$G=e+\g g.$$
From now on we shall suppose that the ground field $k$ is of characteristic zero. For any $x\in\Cal H^n$ the exponential~:
$$e^{*\alpha}(x)=\sum_{k\ge 0}{\alpha^{*k}(x)\over k!}$$
is a finite sum (ending up at $k=n$). It is a bijection from $\g g$ onto $G$. Its inverse is given by~:
$$\mop{Log}(1+\alpha)(x)=\sum_{k\ge 1}{(-1)^{k-1}\over k}\alpha^{*k}(x).$$
This sum again ends up at $k=n$ for any $x\in\Cal H^n$. Let us introduce a decreasing filtration on $\Cal L=\Cal L(\Cal H,\Cal A)$~:
$$\Cal L^n=\{\alpha\in\Cal L, \alpha\restr{\Cal H^{n-1}}=0\}.$$
Clearly $\Cal L_0=\Cal L$ and $\Cal L_1=\g g$. We define the valuation $\mop{val}\varphi$ of an element $\varphi$ of $\Cal L$ as the biggest integer $k$ such that $\varphi$ is in $\Cal L_k$. We shall consider in the sequel the ultrametric distance on $\Cal L$ induced by the filtration~:
$$d(\varphi,\psi)=2^{-\smop{val}(\varphi-\psi)}.$$
For any $\alpha,\beta\in\g g$ let $[\alpha,\beta]=\alpha*\beta-\beta*\alpha$.
\prop{II.3.3}
1) We have the inclusion~:
$$\Cal L^p*\Cal L^q\subset\Cal L^{p+q}.$$
2) The metric space $\Cal L$ endowed with then distance defined just above is complete.
\dem
Take any $x\in\Cal H^{p+q-1}$, and any $\alpha\in\Cal L_p$ and $\beta\in\Cal L_q$. We have~:
$$(\alpha*\beta)(x)=\sum_{(x)}\alpha(x_1)\beta(x_2).$$
Recall that we denote by $|x|$ the minimal $n$ such that $x\in\Cal H^n$. Since $|x_1|+|x_2|=|x|\le p+q-1$, either $|x_1|\le p-1$ or $|x_2|\le q-1$, so the expression vanishes. Now if $(\psi_n)$ is a Cauchy sequence in $\Cal L$ it is immediate to see that this sequence is {\sl locally stationary\/}, i.e. for any $x\in\Cal H$ there exists $N(x)\in\N$ such that $\psi_n(x)=\psi_{N(x)}(x)$ for any $n\ge N(x)$. Then the limit of $(\psi_n)$ exists and is clearly defined by~:
$$\psi(x)=\psi_{N(x)}(x).$$  
\qed
As a corollary the Lie algebra $\Cal L_1=\g g$ is {\sl pro-nilpotent}, in a sense that it is the projective limit of the Lie algebras $\g g/\Cal L^n$, which are nilpotent. 
\alinea{II.4. Algebra morphisms and cocycles}
Let $\Cal H$ be a connected filtered Hopf algebra over $k$, and let $\Cal A$ be a $k$-algebra. A {\sl cocycle from $\Cal H$ to $\Cal A$\/} is a linear morphism $\tau:\Cal H\to\Cal A$ such that $\tau (xy)=\tau(yx)$ for any $x,y\in\Cal H$. It is indeed a $1$-cocycle in the cohomology of the Lie algebra $\Cal H$ with values in $\Cal A$ considered as a trivial $\Cal H$-module. In the case where $\Cal A$ is the ground field $k$ cocycles are just traces.
\ssq
We shall also consider algebra morphisms from $\Cal H$ to $\Cal A$. When the algebra $\Cal A$ is commutative we shall call them slightly abusively {\sl characters\/}. It is clear that any character in our sense is a cocycle. We recover of course the usual notion of character when the algebra $\Cal A$ is the ground field $k$.
\ssq
The notions of character and cocycle involve only the algebra structure of $\Cal H$. Let us consider now the full Hopf algebra structure and see what happens with the convolution product~:
\prop{II.4.1}
Let $\Cal H$ be a connected filtered Hopf algebra over $k$, and let $\Cal A$ be a $k$-algebra. Then,
\medskip
1). The convolution of two cocycles in $\Cal L(\Cal H,\Cal A)$ is a cocycle.
\medskip
2). If $\tau$ is a cocycle such that $\tau(\un)=\un_{\Cal A}$, then the inverse $\tau^{*-1}$ is a cocycle as well.
\medskip
3). In the case of a commutative algebra $\Cal A$ the characters from $\Cal H$ to $\Cal A$ form a group $\Gamma$ under the convolution product, and for any $\chi\in \Gamma$ the inverse is given by~:
$$\chi^{*-1}=\chi\circ S.$$
\dem
Using the fact that $\Delta$ is an algebra morphism we have for any $x,y\in\Cal H$~:
$$f*g(xy)=\sum_{(x)(y)}f(x_1y_1)g(x_2y_2).$$
If $f$ and $g$ are cocycles we get~:
$$\eqalign{f*g(xy)	&=\sum_{(x)(y)}f(y_1x_1)g(y_2x_2)	\cr
			&=f*g(yx).	\cr}$$
If $\Cal A$ is commutative and if $f$ and $g$ are characters we get~:
$$\eqalign{f*g(xy)	&=\sum_{(x)(y)}f(x_1)f(y_1)g(x_2)g(y_2)	\cr
			&=\sum_{(x)(y)}f(x_1)g(x_2)f(y_1)g(y_2)	\cr
			&=(f*g)(x)(f*g)(y).	\cr}$$
The unit $e=u_{\Cal A}\varepsilon$ is both a cocycle and an algebra morphism. The formula for the inverse of a character comes easily from the commutativity of the following diagram~:
\diagramme{
\xymatrix{&\Cal H\otimes\Cal H	\ar[rr]^{S\otimes I}
				&&\Cal H\otimes \Cal H\ar[dr]^{m}\ar[rr]^{f\otimes f}	&& \Cal\otimes\Cal A	\ar[dr]^{m_{\Cal A}}	&\\
\Cal H\ar[rr]_\varepsilon \ar[dr]^\Delta \ar[ur]^\Delta
	\ar@/_4.5pc/	@{-->}[rrrrrr]_{f*(f\circ S)}
	\ar@/^4.5pc/	@{-->}[rrrrrr]^{(f\circ S)*f}
	\ar@/^1pc/	@{-->}[rrrrrr]_e
				&&	k\ar[rr]_u	&&\Cal H\ar[rr]_f
	&&	\Cal A	\\
&\Cal H\otimes\Cal H	\ar[rr]^{I\otimes S}
				&&\Cal H\otimes \Cal H\ar[ur]^{m}
\ar[rr]^{f\otimes f}	&&\Cal A\otimes\Cal A\ar[ur]^{m_\Cal A}	&}
}
Finally the fact that the inverse of a cocycle $\tau$ such that $\tau(\un)=\un_{\Cal A}$ is a cocycle comes from 1) and from the formula~:
$$\tau\inver(x)=\sum_{k\ge 0}(e-\tau)^{*k}(x).$$
\qed
We call {\sl derivations \/ \rm(or \sl infinitesimal characters\/\rm) \sl with values in the algebra $\Cal A$} those elements $\alpha$ of $\cal L(\Cal H,\Cal A)$ such that~:
$$\alpha(xy)=e(x)\alpha(y)+\alpha(x)e(y).$$

\prop{II.4.2}
Suppose that $\Cal A$ is a commutative algebra.
Let $G_1$ (resp. $\g g_1$) be the set of characters of $\Cal H$ with values in $\Cal A$ (resp the set of derivations of $\Cal H$ with values in $\Cal A$), and let $G_2$ (resp. $\g g_2$) be the set of cocycles $\varphi$ from $\Cal H$ to $\Cal A$ such that $\varphi(\un)=\un_{\Cal A}$ (resp. $\varphi(\un)=0$). Then $G_1$ and $G_2$ are subgroups of $G$, the exponential restricts to a bijection from $\g g_1$ onto $G_1$ (resp. from $\g g_2$ onto $G_2$), and $\g g_1,\g g_2$ are Lie subalgebras of $\g g$.
\dem
Part of these results are a reformulation of proposition II.4.1 and some points are straightforward. The only non-trivial point concerns $\g g_1$ and $G_1$. Take two derivations $\alpha$ and $\beta$ with values in $\Cal A$ and compute~:
$$\eqalign{(\alpha*\beta)(xy)
	&=\sum_{(x)(y)}\alpha(x_1x_2)\beta(y_1y_2)	\cr
	&=\sum_{(x)(y)}\bigl(\alpha(x_1)e(y_1)+e(x_1)\alpha(y_1)\bigr).
		\bigl(\beta(x_2)e(y_2)+e(x_2)\alpha(y_2)\bigr)		\cr
	&=(\alpha*\beta)(x)e(y)+\alpha(x)\beta(y)+\beta(x)\alpha(y)
		+e(x)(\alpha*\beta)(y).	\cr}$$
Using the commutativity of $\Cal A$ we immediately get~:
$$[\alpha,\beta](xy)=[\alpha,\beta](x)e(y)+e(x)[\alpha,\beta](y),$$
which shows that $\g g_1$ is a Lie algebra. Now for $\alpha\in\g g_1$ we have~:
$$\alpha^{*n}(xy)=\sum_{k=0}^n{n\choose k}\alpha^{*k}(x)\alpha^{*(n-k)}(y),$$
as easily see by induction on $n$. A straightforward computation then yields~:
$$e^{*\alpha}(xy)=e^{*\alpha}(x)e^{*\alpha}(y).$$
\qed

\alinea{II.5. Birkhoff decomposition}
We consider here the situation where the algebra $\Cal A$ admits a {\sl
renormalization scheme\/}, i.e. a splitting into two subalgebras~:
$$\Cal A=\Cal A_-\oplus \Cal A_+$$
with $\un\in\Cal A_+$. As an example, take $\Cal A$ as the field $k[t\inver,t]]$ of Laurent series, $\Cal A_-=t\inver k[t\inver]$ and $\Cal A_+=k[[t]]$. The projection on $\Cal A_-$ parallel to $\Cal A_+$ will be denoted by $\pi$.
\th{II.5.1}
1). Let $\Cal H$ be a connected filtered Hopf algebra. Let $G$ be the group of those $\varphi\in\Cal L(\Cal H,\Cal A)$ such that $\varphi(\un)=\un_{\Cal A}$ endowed with the convolution product. Any $\varphi\in G$ admits a unique Birkhoff decomposition~:
$$\varphi=\varphi_-^{*-1} * \varphi_+,$$
where $\varphi_-$ sends $\un$ to $\un_{\Cal A}$ and $\mop{Ker}\varepsilon$ into $\Cal A_-$, and where $\varphi_+$ sends $\Cal H$ into $A_+$. The maps $\varphi_-$ and $\varphi_+$ are given on $\mop{Ker}\varepsilon$ by the following recursive formulas~:
$$\eqalign{\varphi_-(x)	&=-\pi\Bigl( \varphi(x)+\sum_{(x)}\varphi_-(x')\varphi(x'')\Bigr)	\cr
\varphi_+(x)	&=(I-\pi)\Bigl( \varphi(x)+\sum_{(x)}\varphi_-(x')\varphi(x'')\Bigr).	\cr}$$
2). If $\tau\in G$ is a cocycle, the components $\tau_-$ and $\tau_+$ occurring in the Birkhoff decomposition of $\tau$ are cocycles as well.
\bigskip
3). If the algebra $\Cal A$ is commutative and if $\chi$ is a character, the components $\chi_-$ and $\chi_+$ occurring in the Birkhoff decomposition of $\chi$ are characters as well. 
\dem
Points 1) and 3) together give an abstract counterpart of Theorem 4 of \cite {CK}, point 2) is new up to my knowledge. The proof goes along the same lines~: for the first assertion it is immediate from the definition of $\pi$ that $\varphi_-$ sends $\mop{Ker}\varepsilon$ into $\Cal A_-$, and that $\varphi_+$ sends $\mop{Ker}\varepsilon$ into $\Cal A_+$. It only remains to check equality $\varphi_+=\varphi_-*\varphi$, which is an easy computation~:
$$\eqalign{\varphi_+(x)	&=(I-\pi)\Bigl( \varphi(x)+\sum_{(x)}\varphi_-(x')\varphi(x'')\Bigr).	\cr
			&=\varphi(x)+\varphi_-(x)+ \sum_{(x)}\varphi_-(x')\varphi(x'')	\cr
			&=(\varphi_-*\varphi)(x).\cr}$$
To prove the second assertion it is sufficient to prove that $\tau_-$ is a cocycle whenever $\tau$ is a cocycle. The same property for $\tau_+$ comes then from proposition II.3.1. We prove the formula $\tau_-(xy)=\tau_-(yx)$ by induction on the integer $d=|x|+|y|$~: it is true for $d\le 1$. Suppose the formula is true up to $d-1$ and take any $x,y\in\Cal H$ with $|x|+|y|=d$. Decompose $\Delta(xy)$ with the second version of Sweedler's notation~:
$$\displaylines{\Delta(xy)=xy\otimes\un+\un\otimes xy+x\otimes y+y\otimes x\hfill\cr
\hfill +\sum_{(x)}(x'y\otimes x''+x'\otimes x''y)
	+\sum_{(y)}(xy'\otimes y''+y'\otimes xy'')\hfill	\cr
\hfill +\sum_{(x)(y)}x'y'\otimes x''y''.	\cr}$$
We have then~:
$$\displaylines{\tau_-(xy)=-\pi\Bigl(\tau(xy)+\tau_-(x)\tau(y)+\tau_-(y)\tau(x)\hfill\cr
\hfill
+\sum_{(x)}\bigl(\tau_-(x'y)\tau(x'')+\tau_-(x')\tau(x''y) \bigr)
+\sum_{(y)}\bigl(\tau_-(xy')\tau(y'')+\tau_-(y')\tau(xy'') \bigr) \hfill\cr
\hfill +\sum_{(x)(y)}\tau_-(x'y')\tau(x''y'')	\Bigr),\cr}$$
whereas~:
$$\displaylines{\tau_-(yx)=-\pi\Bigl(\tau(yx)+\tau_-(y)\tau(x)+\tau_-(x)\tau(y)\hfill\cr
\hfill
+\sum_{(y)}\bigl(\tau_-(y'x)\tau(y'')+\tau_-(y')\tau(y''x) \bigr)
+\sum_{(x)}\bigl(\tau_-(yx')\tau(x'')+\tau_-(x')\tau(yx'') \bigr) \hfill\cr
\hfill +\sum_{(x)(y)}\tau_-(y'x')\tau(y''x'')	\Bigr).\cr}$$
Using the cocycle property for $\tau$ and the induction hypothesis we see that the two expressions are the same.
\msq
The proof of assertion 3) goes exactly as in \cite {CK} and relies on the following {\sl Rota-Baxter\/} equality in $\Cal A$~:
$$\pi(a)\pi(b)=-\pi(ab)+\pi\bigl(\pi(a)b\bigr)+\pi\bigl(\pi(b)a\bigr),$$
which is easily verified by decomposing $a$ and $b$ into their $\Cal A_\pm$-parts. Let $\chi$ be a character of $\Cal H$ with values in $\Cal A$. Suppose that we have $\chi_-(xy)=\chi_-(x)\chi_-(y)$ for any $x,y\in\Cal H$ such that $|x|+|y|\le d-1$, and compute for $x,y$ such that $|x|+|y|=d$~:
$$\chi_-(x)\chi_-(y)=\pi(X)\pi(Y),$$
with $X=\chi(x)-\sum_{(x)}\chi_-(x')\chi(x'')$ and $Y=\chi(y)-\sum_{(y)}\chi_-(y')\chi(y'')$. Using the formula~:
$$\pi(X)=-\chi_-(x),$$
we get~:
$$\chi_-(x)\chi_-(y)=-\pi\bigl(XY+\chi_-(x)Y+X\chi_-(y)\bigr),$$
hence~:
$$\displaylines{\chi_-(x)\chi_-(y)=-\pi\Bigl(\chi(x)\chi(y)+\chi_-(x)\chi(y)+\chi(x)\chi_-(y)\hfill\cr
\hfill +\sum_{(x)}\chi_-(x')\chi(x'')\bigl(\chi(y)+\chi_-(y)\bigr)
	+\sum_{(y)}\bigl(\chi(x)+\chi_-(x)\bigr)\chi_-(y')\chi(y'')\hfill\cr
\hfill +\sum_{(x)(y)}\chi_-(x')\chi(x'')\chi_-(y')\chi(y'')\Bigl).\cr}$$
We have to compare this expression with~:
$$\displaylines{\chi_-(xy)=-\pi\Bigl(\chi(xy)+\chi_-(x)\chi(y)+\chi_-(y)\chi(x)\hfill\cr
\hfill
+\sum_{(x)}\bigl(\chi_-(x'y)\chi(x'')+\chi_-(x')\chi(x''y) \bigr)
+\sum_{(y)}\bigl(\chi_-(xy')\chi(y'')+\chi_-(y')\chi(xy'') \bigr) \hfill\cr
\hfill +\sum_{(x)(y)}\chi_-(x'y')\chi(x''y'')	\Bigr).\cr}$$
These two expressions are easily seen to be equal using the commutativity of the algebra $\Cal A$, the character property for $\chi$ and the induction hypothesis.
\qed
{\bf Remark} : define the {\sl Bogoliubov character\/} as the map $b:G\to \Cal
L(\Cal H,\Cal A)$ recursively given by~:
$$b(\varphi)(x)=\varphi(x)+\sum_{(x)}\varphi_-(x')\varphi(x'').$$
Then the components of $\varphi$ in the birkhoff decomposition read~:
$$\varphi_-=-\pi\circ b(\varphi),\hskip 12mm \varphi_+=(I-\pi)\circ b(\varphi).$$
\alinea{II.6. The BCH approach to Birkhoff decomposition}
\qquad Let $\g g$ a Lie algebra endowed with a decreasing filtration~:
$$\g g=\g g_1\supset\g g_2\supset\cdots\supset \g g_n\supset\cdots$$
such that the intersection of the $\g g_i$'s is reduced to $\{0\}$. We ask for the inclusion~:
$$[\g g_i,\g g_j]\subset\g g_{i+j},$$
so in particular $\g g$ is a pro-nilpotent Lie algebra. Suppose that $\g g$ is complete for the (metric) topology defined by this filtration. The Baker-Campbell-Hausdorff series defines then a pro-nilpotent group law on $\g g$~:
$$X.Y=X+Y+\frac 12[X,Y]+\frac 1{12}([X,[X,Y]]+[Y,[Y,X]])+\cdots$$
The following proposition is implicitly used in \cite {EGK2}~:
\prop{II.6.1}
For any linear map $R:\g g\to\g g$ preserving the filtration there exists a (usually non-linear) map $\chi_R:\g g\to\g g$ such that $(\chi_R-\mop{Id}_{\sg g})(\g g_i)\subset \g g_{2i}$ for any $i\ge 1$, and such that, with $\wt R:=\mop{Id}_{\sg g}-R$ we have~:
$$\forall X\in\g g,\ \ X=R\bigl(\chi_R(X)\bigr).\wt R(\chi_R(X)\bigr).
\eqno{(*)}$$
\dem
Let us introduce for any $X,Y\in\g g$ the following expression~:
$$\delta(X,Y)=X.Y-X-Y
	=\frac 12[X,Y]+\frac 1{12}([X,[X,Y]]+[Y,[Y,X]])+\cdots$$
Then equation (*) can be rewritten as~:
$$\chi_R(X)=F_X\bigl(\chi_R(X)\bigr),$$
with $F_X:\g g\to\g g$ defined by~:
$$F_X(Y)=X-\delta\bigl(R(Y),\wt R(Y)\bigr).$$
This map $F_X$ is a contraction with respect to the metric associated with the filtration~: indeed if $Y,\varepsilon\in\g g$ with $\varepsilon\in\g g_n$, we have~:
$$F_X(Y+\varepsilon)-F_X(Y)=\delta\bigl(R(Y),\wt R(Y)\bigr)
	-\delta\bigl(R(Y+\varepsilon),\wt R(Y+\varepsilon)\bigr).$$
Right-hand side is a sum of iterated commutators in each of which $\varepsilon$ does appear at least once. So it belongs to $\g g_{n+1}$. So the sequence $F_X^n(Y)$
converges in $\g g$ to a unique fixed point $\chi_R(X)$ for $F_X$.
\ssq
Let us remark that for any $X\in\g g_i$, the element $F_X(X)-X$ belongs to $\g g_{2i}$. Now taking $X$ as starting point it is obvious from the expression~:
$$\chi_R(X)-X=\sum_{k=1}^{+\infty}(F_X^k-F_X^{k-1})(X)$$
that for any $X\in\g g_i$, the element $\chi_R(X)-X$ belongs to $\g g_{2i}$.
\qed
Let $\Cal H$ be a connected filtered Hopf algebra, let $\Cal A$ a commutative algebra endowed with a splitting $\Cal A=\Cal A_-\oplus\Cal A_+$ as in \S\ II.5, and let $\pi$ be the projection on $\Cal A_-$ parallel to $\Cal A_+$. Now take for $\g g$ any of the Lie algebras $\g g,\g g_1,\g g_2$ of \S\ II.3 and II.4, set $G=\exp_*(\g g)$ and set $R(X)=\pi\circ X$. Of course $R$ makes sense on $\Cal L(\Cal H,\Cal A)$. Equation (*) then yields the following equality in the group $G$, $G_1$ or $G_2$~:
$$e^{*X}=e^{*R\bigl(\chi_R(X)\bigr)}*e^{*\wt R\bigl(\chi_R(X)\bigr)},\eqno{(**)}$$
as map $R$ respects the decreasing filtration introduced in \S\ II.3. Now, thanks to Rota-Baxter relation in $\Cal A$~:
$$\pi(a)\pi(b)=\pi\Bigl(\pi(a)b+\pi(b)a-ab\Bigr),$$
K. Ebrahimi-Fard, L. Guo and D. Kreimer derive in \cite {EGK2} two identities
involving the Bogoliubov character~:
$$e^{*-R\bigl(\chi_R(X)\bigr)}=R\bigl(b(e^{*X})\bigr),\hskip 12mm e^{*\wt R\bigl(\chi_R(X)\bigr)}=-\wt R\bigl(b(e^{*X})\bigr).$$
So (**) is indeed the Birkhoff decomposition of the element $\varphi=e^{*X}$ of $G$, $G_1$ or $G_2$, namely~:
$$\varphi_-=e^{*-R\bigl(\chi_R(X)\bigr)},\hskip 12mm 
	\varphi_+=e^{*\wt R\bigl(\chi_R(X)\bigr)}.$$
Of course when $\Cal H$ is cocommutative the convolution product is commutative, the Lie algebras involved are abelian and the situation simplifies greatly, as $\chi_R=\mop{Id}_{\sg g}$ for any linear map $R$ here \cite {EGK1}. Further developments on general Rota-Baxter algebras can be found in \cite {EGK2} and \cite {EGK3}.
\smallskip
{\bf Remark}~: Rota-Baxter identity for $R$ just guarantees that equation (**) gives a Birkhoff decomposition. When $R$ is idempotent (which is indeed the case in the setting of \S\ II.5, where $R$ is built up from projection $\pi$) this decomposition is unique and is given either by (**) or by the recursive formulas of Theorem II.5.1. On the other hand, with an idempotent $R$ which does {\sl not\/} verify Rota-Baxter identity (or, which is the same, with a direct sum decomposition of the algebra $\Cal A$ on two components $\Cal A_+$ and $\Cal A_-$ which are {\sl not\/} subalgebras of $\Cal A$), we still get a unique Birkhoff decomposition along the lines of theorem II.5.1 for the groups $G$ and $G_2$, which does not coincide with (**). But Rota-Baxter identity arises in an essential way, as we have seen, in order to get the Birkhoff decomposition for the groups $G_1$ of $\Cal A$-valued characters.
\alinea{II.7. Renormalized traces and characters}
\qquad Keeping the same notations we take $k=\C$ as ground field and we specialize to the case when $\Cal A$ is the field of meromorphic functions (resp. the field of germs of meromorphic functions at $z_0$), $\Cal A_+$ is the algebra of meromorphic functions which are holomorphic at $z_0$ (resp. the algebra of germs of holomorphic functions at $z_0$) and $\Cal A_-$ is the non-unital algebra $(z-z_0)\inver \C[(z-z_0)\inver]$. Applying the projection $\pi$ to such a meromorphic fuction amounts to ``take its divergent part at $z_0$''. This particular splitting is called the {\sl minimal subtraction scheme\/}. It is by no means unique~: for example for any automorphism $\theta$ of the field of (germs of) meromorphic functions, we can consider the splitting~:
$$\Cal A=\Cal A^\theta_+\oplus\Cal A^\theta_-,$$
with $\Cal A_\pm^\theta=\theta(\Cal A_\pm)$. The corresponding projection on $\Cal A_-^\theta$ is given by $\pi^\theta=\theta\circ\pi\circ\theta\inver$. As an example of automorphism $\theta$ fix a constant $c$, consider the change of variable $z\mapsto z'$ such that~:
$${1\over z'-z_0}={1\over z-z_0}+c$$
(hence $z'=z_0+\displaystyle{z-z_0\over 1+c(z-z_0)}$), and set $\theta(f)(z)=f(z')$.
\msq
Considering a linear map $\varphi$ from $\Cal H$ to $\Cal A$ (and a particular splitting of $\Cal A$) we can consider its Birkhoff decomposition $\varphi=\varphi_-^{*-1}*\varphi_+$ given by Theorem II.5.1, and evaluate $\varphi_+(x)$ at $z=z_0$ for any $x\in\Cal H$. This gives a linear map $\varphi_+^{z_0}$ from $\Cal H$ to $\C$ which we call the {\sl renormalized value of $\varphi$ at $z=z_0$\/}. According to Theorem II.5.1 the renormalized value of a cocycle at $z_0$ is a trace, and the renormalized value of a character at $z=z_0$ is a $\C$-valued character. All this procedure depends in an essential way on the choice of the renormalization scheme, i.e. the splitting of $\Cal A$.
\alinea{II.8. More on connected graded Hopf algebras}
Let $\Cal H$ be a connected graded Hopf algebra. The grading induces a biderivation $Y$ defined on homogeneous elements by~:
$$\eqalign{Y:\Cal H_n	&\longrightarrow \Cal H_n	\cr
		x	&\longmapsto	nx.\cr}$$
Exponentiating we get a one-parameter group $\theta_t$ of automorphisms of the Hopf algebra $\Cal H$, defined on $\Cal H_n$ by~:
$$\theta_t(x)=e^{nt}x.$$
\lemme{II.8.1}
$\varphi\mapsto \varphi\circ Y$ is a derivation of $\bigl(\cal L(\Cal H,\Cal A),*\bigr)$, and $\varphi\mapsto \varphi\circ \theta_t$ is an automorphism of $\bigl(\cal L(\Cal H,\Cal A),*\bigr)$ for any complex $t$. 
\dem
We compute for $\varphi,\psi\in \cal L(\Cal H,\Cal A)$ and for an homogeneous element $x$ of $\Cal H$~:
$$\eqalign{(\varphi*\psi)\circ Y(x)	&=|x|\sum_{(x)}\varphi(x_1)\psi(x_2)\cr
		&=\sum_{(x)}(|x_1|+|x_2|)\varphi(x_1)\psi(x_2)	\cr
		&=\sum_{(x)}(\varphi\circ Y)(x_1)\psi(x_2)
			+\varphi(x_1)(\psi\circ Y)(x_2)		\cr
		&=\bigl((\varphi\circ Y)*\psi + \varphi*(\psi\circ Y)\bigr)(x),\cr
				}$$
which shows the first assertion. The second part of the proposition is proven similarly.
\qed
Using the fact that $e\circ Y=0$ we easily compute for any derivation $\alpha$ with values in $\Cal A$~:
$$\eqalign{(\alpha\circ Y)(xy)	&=\alpha(Y(x).y+x.Y(y))	\cr
	&=(\alpha\circ Y)(x)e(y)+(e\circ Y)(x)\alpha(y)
		+\alpha(x)(e\circ Y)(y)+e(x)(\alpha\circ Y)(y)\cr
	&=(\alpha\circ Y)(x)e(y)+e(x)(\alpha\circ Y)(y).\cr}$$
So we have proved~:
\lemme{II.8.2}
The map $\alpha\mapsto \alpha\circ Y$ is a linear automorphism of the space of derivations of $\Cal H$ with values in $\Cal A$. Its inverse is given by $\alpha\mapsto\alpha\circ Y\inver$, where $Y\inver(x)=|x|\inver x$ for $x$ homogeneous of positive degree, and $Y\inver(\un)=0$.
\ndem
{\bf Remark}~: the notation $Y\inver$ is of course slightly incorrect, as the inverse of $Y$ does not make sense on $\Cal H_0$. The convention $Y\inver(\un)=0$ is arbitrary~: any other value of $Y\inver(\un)$ would give the same result, as derivations with values in $\Cal A$ vanish at $\un$. 
\alinea{II.9. Examples}
{\sl II.9.1. The Hopf algebra of positive integers\/}
\medskip
This example is a simplified version of the one given by D. Kreimer in \cite {K2 \S\ 2.1}. Consider the algebra $\Cal N$ of the multiplicative semigroup $\N^*=\{1,2,3,\ldots\}$ of positive integers. As a vector space it admits a basis $(e_n)_{n\in\N^*}$ with product given by $e_n.e_m=e_{nm}$ and extended by linearity. We endow $\Cal N$ with a structure of commutative cocommutative connected graded Hopf algebra thanks to the decomposition of any integer into a product of prime factors~: namely we set $\Delta(e_1)=e_1\otimes e_1$, and for any prime $p$~:
$$\Delta(e_p)=e_p\otimes e_1+e_1\otimes e_p,$$
and we extend $\Delta$ to an algebra isomorphism. Hence,
$$\Delta(e_{p_1\cdots p_k})=\sum_{I\amalg J=\{1,\ldots ,k\}}e_{p_I}\otimes e_{p_J},$$
where $p_I$ denotes the product of the primes $p_j,j\in I$. The grading is clearly given by the number of prime factors (including multiplicities). The antipode is given by~:
$$S(e_n)=(-1)^{|n|}e_n.$$
Suppose that the ground field is $k=\C$. The map $n\mapsto n^z$ defines a character $\varphi$ of $\Cal N$ with values into the holomorphic functions. Then the Riemann Zeta function is nothing but the evaluation of $\varphi$ on the element~:
$$\omega=e_1+e_2+e_3+\cdots=\prod_{p\hbox{ \sevenrm prime }}
	{1\over 1-e_p}.$$
Here $1/(1-e_p)$ stands for the infinite sum~: $e_1+e_p+e_{p^2}+\cdots$.
Of course $\omega$ is not an element of $\Cal N$ : it makes sense (as well as the abstract Euler product expansion on the right-hand side) only in the completion of $\Cal N$ with respect to the {\sl fine filtration\/} defined by the vector space grading $d(n)=n-1$. But evaluating the character $\varphi$ on both sides of this equality gives the well-known Euler product expression of the Zeta function.
\medskip
{\sl II.9.2. Tensor and symmetric algebras\/}
\medskip
The tensor Hopf algebra $T(V)$ of any vector space $V$ (cf. Example I.6.2) is obviously graded. The symmetric Hopf algebra is a particular case of enveloping Hopf algebra, with $V$ viewed as an abelian Lie algebra. The Hopf algebra $S(V)$ is a cocommutative commutative connected graded Hopf algebra. Note that an enveloping algebra is not graded in general, since the quotienting ideal generated by $x\otimes y-y\otimes x-[x,y]$ is not homogeneous.
\medskip
{\sl II.9.3. Planar decorated rooted trees\/}
\medskip
We borrow in this section some material from \cite {F}. A {\sl planar rooted tree\/} is an oriented connected contractible graph, with a finite number of vertices, together with an embedding of it into the plane, such that only one vertex has only outgoing edges (the root). We have drawn below the planar rooted trees with four vertices~:
\dessin{30mm}{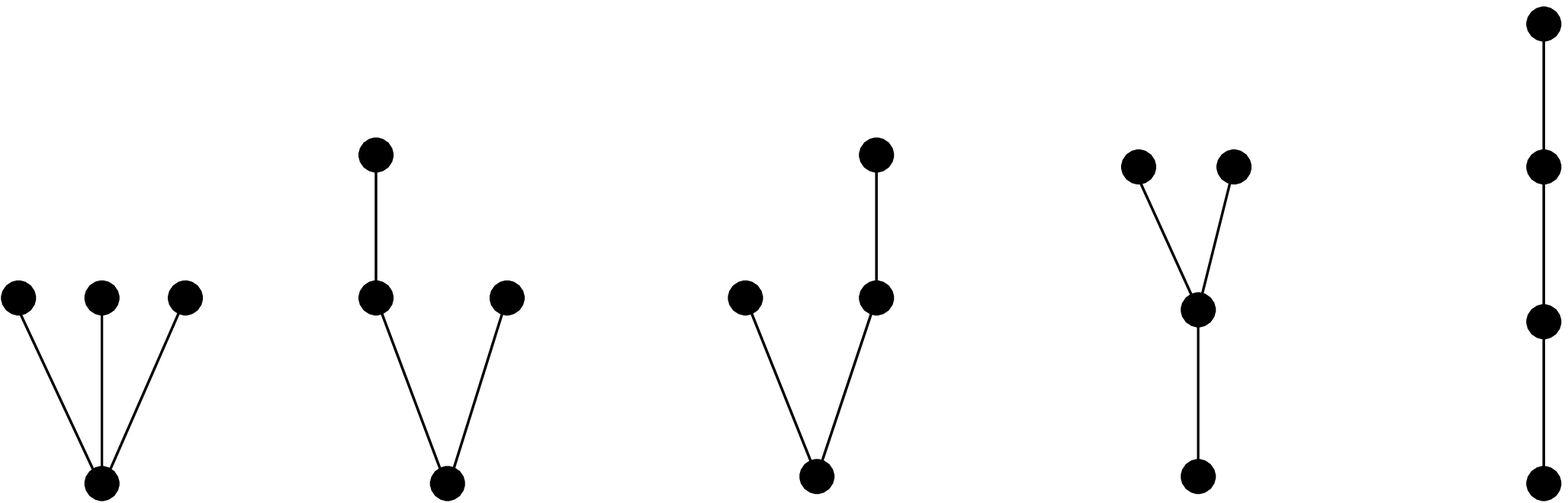}{1}{the planar rooted trees with four vertices}
Let $\Cal T$ be the set of planar rooted trees. Let $V$ be a vector space on some field $k$, and let $t$ be a planar rooted tree. The {\sl space of decorations of $t$ by $V$\/} is the vector space $V^{\otimes t}$. A planar rooted tree $t$ together with an indecomposable element of $V^{\otimes t}$ is called a {\sl decorated rooted tree\/}. Let us consider the vector space~:
$$\Cal T_V=\bigoplus_{t\in\Cal T}V^{\otimes t},$$
let $\Cal H_V$ be the (noncommutative) free algebra generated by $\Cal T_V$. Products of decorated trees (decorated forests) generate $\Cal H_V$ as a graded vector space, the degree of a decorated forest being given by the total number of vertices. The connected graded Hopf algebra structure on $\Cal H_V$ is given by the co-unit $\varepsilon$ sending $\un$ to $1$ and any nonempty decorated forest to $0$, and by a coproduct which we describe shortly here~:
\ssq
An {\sl elementary cut\/} on a tree is a cut on some edge of the given tree. An {admissible cut\/} is a cut such that any path starting from the root contains at most one elementary cut. The {\sl empty cut\/} is considered as elementary, as well as the {\sl total cut\/}, i.e. a cut below the root.
A cut on a forest is said to be admissible if its restriction to any tree factor is admissible. Any elementary cut $c$ sends a forest $F$ to a couple $\bigl(P^c(F),R^c(F)\bigr)$, the {\sl crown\/} and the {\sl trunk\/} respectively. The trunk of a tree is a tree, but the crown of a tree is a forest. Let $\mop{Adm}(F)$ the set of admissible cuts of the forest $F$, and let $\mop{Adm}^*(F)$ the set of elementary cuts discarding the empty cut and the total cut. The coproduct~:
$$\Delta(F)=\sum_{c\in\smop{Adm}F} P^c(F)\otimes R^c(F)$$
is graded, co-associative and compatible with the product \cite{F}. The compatibility with the product is clear (due to the definition of an admissible cut for a forest). There is a beautiful proof of the co-associativity in \cite {F} using induction on the degree and grafting of any forest on a decorated root. We propose here a more intuitive proof~: say that a couple $(c_1,c_2)$ of cuts is {\sl bi-admissible\/} if both cuts $c_1,c_2$ are admissible and if $c_1$ never bypasses $c_2$, i.e. if $c_2$ never cuts the trunk of $c_1$. Any bi-admissible couple $c=(c_1,c_2)$ of cuts $c$ defines a crown $P^c(F)=P^{c_2}(F)$, a trunk $R^c(F)=R^{c_1}(F)$, and a middle $M^c(F)$~:
\dessin{25mm}{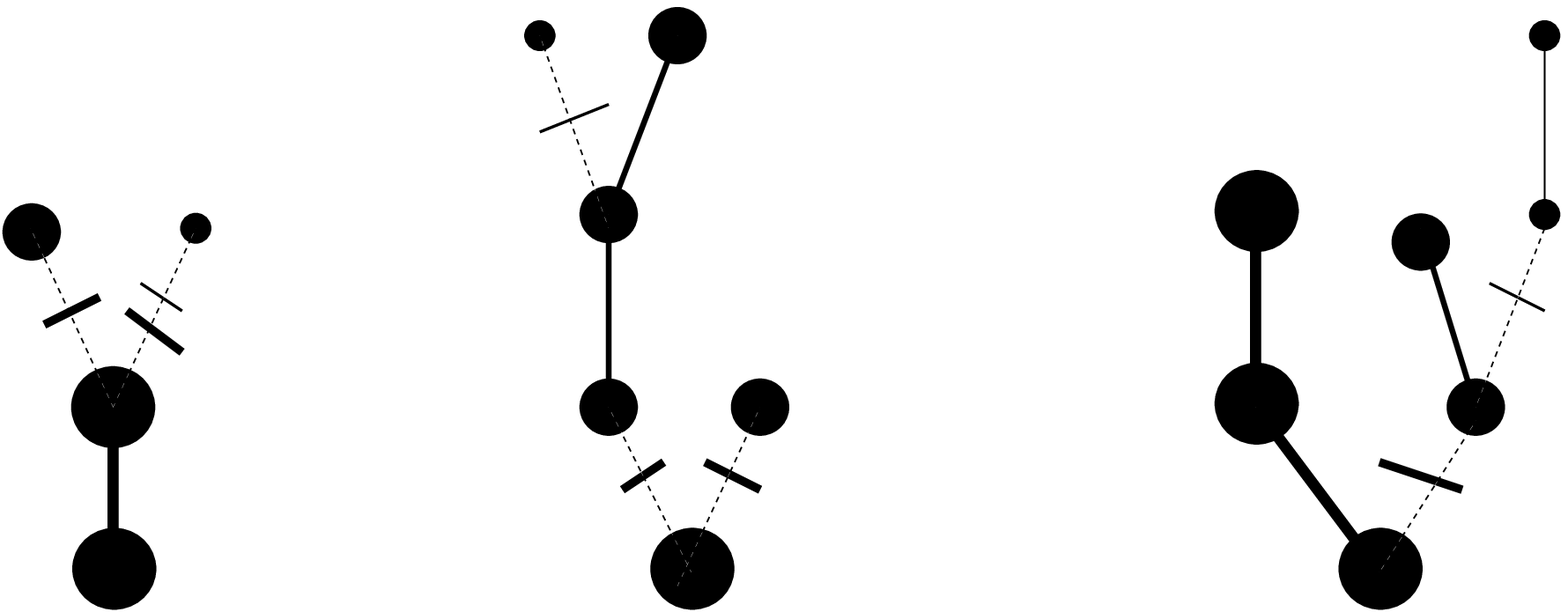}{2}{an example of bi-admissible couple of cuts. From thickest to thinnest : trunk, middle and crown}
 Let $\mop{Adm}_2F$ the set of bi-admissible couples of cuts of the forest $F$. It is quite straightforward to set down the formula for the iterated coproduct~:
$$(\Delta\otimes I)\circ\Delta(F)=(I\otimes\Delta)\circ\Delta(F)
=\sum_{c\in\smop{Adm}_2F} P^c(F)\otimes M^c(F)\otimes R^c(F).$$
Of course the $n$-fold iterated coproduct admits a similar expression, involving $n$-admissible n-uples of admissible cuts and $n+1$ ``level segments'' of the forest, from the crown down to the trunk.  
\ssq
By corollary II.2.2 the connected graded bialgebra $\Cal H_V$ thus obtained admits an antipode given on $\mop{Ker}\varepsilon$ by any of the two recursive formulas~:
$$\eqalign{S(F)	&=-F-\sum_{c\in\smop{Adm}^*(F)} S\bigl(P^c(F)\bigr).R^c(F)	\cr
		&=-F-\sum_{c\in\smop{Adm}^*(F)} P^c(F).S\bigl(R^c(F)\bigr).	\cr}$$
The square of the antipode does not in general coincide with the identity.
\medskip
{\sl II.9.3. Decorated rooted trees\/}
\medskip
The construction is the same except that we consider rooted trees independently from any embedding into the plane, and we consider the {\sl free commutative} algebra generated by decorated rooted trees. We thus obtain a commutative Hopf algebra $\Cal H'_V$ which is clearly a quotient of $\Cal H_V$. This Hopf algebra is thoroughly investigated in \cite {F}.
\paragraphe{III. Hopf algebras of Feynman graphs}
\qquad We treat this example (more exactly this family of examples) in a separate section for two main reasons~: firstly the Hopf algebras appearing there are pointed but not connected, and secondly this is the very example where a link is established with quantum field theory. The non-connectedness is not a very serious problem~: as we shall see we can reason on a connected quotient and go back. The formula for the coproduct will differ slightly from that of Connes-Kreimer in order to deal with this non-connectedness problem, but both will agree on the connected quotient. We follow \cite {K1} quite closely, with some modifications in order to allow self-loops.
\alinea{III.1. Discarding exterior structures}
\qquad
Anyone a little bit familiar with quantum field theory knows that Feynman
graphs are made of internal and external edges of different types, and that an
external edge comes with a vector attached to it (an {\sl exterior
momentum}). The sum of all exterior momenta of a given graph must be equal to
zero, reflecting the global conservation of momenta in an interaction. The
{\sl Feynman rules\/} attach to a graph together with such an external
structure an integral which can be divergent. This integral can be regularized
by various procedures, among them {\sl dimensional regularization\/}~: the
idea is to ``let the dimension of the space of momenta vary in the complex
numbers'', a procedure which has been recently given a precise geometrical
contents by A. Connes and M. Marcolli (\cite{CM2} \S\ 15). The divergent integral is now replaced by a meromorphic function with poles at least at the entire dimensions where the original integral diverges \cite{C Chap. 4}, \cite{E}.
\ssq
The approach of renormalization by A. Connes and D. Kreimer can be summarized as follows~: organize Feynman graphs with their exterior structures into a graded Hopf algebra, understand the (regularized, e.g. by means of dimensional regularization) Feynman rules as a character of this Hopf algebra with values into some algebra $\Cal A$ (e.g. the meromorphic functions), choose a renormalization scheme, i.e. a splitting $\Cal A=\Cal A_+\oplus\Cal A_-$ into two subalgebras, apply the method of \S\ II.5 and II.6 to extract a renormalized value, and finally recognize that this method agrees with algorithms already developed by physicists, such as the Bogoliubov-Parasiuk-Hepp-Zimmerman (BPHZ) algorithm.
\msq
Our first step will consist in constructing a Hopf algebra from Feynman diagrams without exterior structure (i.e. with exterior momenta nullified). 
\alinea{III.2. Operations on Feynman graphs} 
\qquad
A {\sl Feynman graph\/} is a (non-oriented, non-planar) graph with a finite number of vertices and edges. An {\sl internal edge\/} is an edge connected at both ends to a vertex (which can be the same in case of a self-loop), an {\sl external edge\/} is an edge with one open end, the other end being connected to a vertex.
\ssq
A Feynman graph is called by physicists {\sl vacuum graph, tadpole graph, self-energy graph\/}, resp. {\sl interaction graph\/} if its number of external edges is $0$,$1$, $2$, resp. $>2$.
\ssq
The edges (internal or external) will be of different types labelled by a positive integer ($1,2,3,\ldots$), each type being represented by the way the corresponding edge is drawn (full, dashed, wavy, various colours, etc...). Let $\tau (e)\in\N^*$ be the type of the edge $e$. For any vertex $v$ let $\mop{st}(v)$ be the {\sl star\/} of $v$, i.e. the set of all edges attached to $v$, {\sl with self-loops counted twice\/}. Hence the valence of the vertex is given by the cardinal of $\mop{st}(v)$. Finally to each vertex we associate its {\sl type\/}, the sequence $(\uple nr)$ of positive integers where $n_j$ stands for the number of edges of type $j$ in $\mop{st}(v)$. Let $T(v)$ be the type of the vertex $v$
\ssq
For example, in the {\sl $\varphi^n$ theory\/} there is only one type of edge, and two types of vertices~: the bivalent vertices and the $n$-valent vertices. In quantum electrodynamics there are two types of edges : the fermion edges (usually drawn full), and the boson edges (usually drawn wavy), and three types of vertices~: bivalent boson-boson vertices, bivalent fermion-fermion vertices, and trivalent vertices with two fermion edges and one boson edge. Most of the pictures will be drawn in $\varphi^3$ or $\varphi^4$ theory, or in quantum electrodynamics. 
\ssq
A {\sl one-particle irreducible graph\/} (in short, 1PI graph) is a connected graph which remains connected when we cut any internal edge. A disconnected graph is said to be {\sl locally 1PI\/} if any of its connected components is 1PI. The {\sl residue\/} of a connected graph is the graph with only one vertex obtained by shrinking all internal edges to a point. Of course any connected graph has the same type as its residue.
\dessin{30mm}{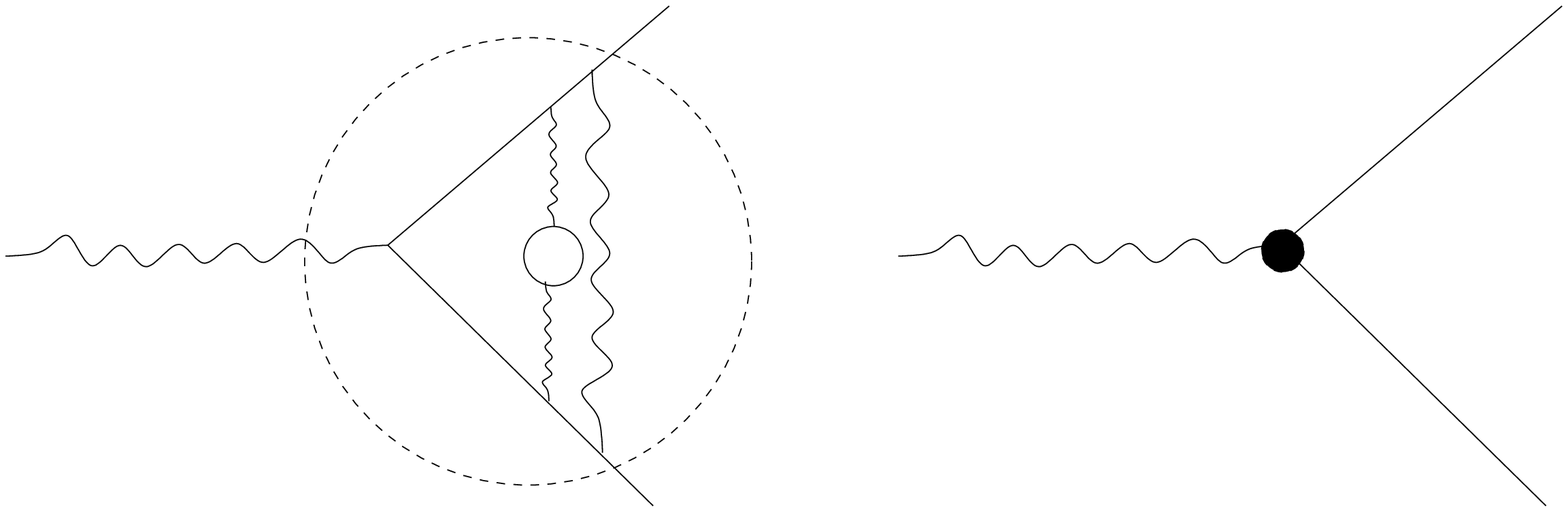}{3}{A QED interaction graph and its residue}
A {\sl subgraph\/} of a Feynman graph is either the empty graph, or a {\sl nonempty\/} (connected or disconnected) set of internal edges together with the vertices they encounter and the stars of those vertices. A {\sl proper subgraph\/} of $\Gamma$ is a subgraph different from the empty graph or the whole graph $\Gamma$ itself. If $\gamma$ is a subgraph inside a graph $\Gamma$, the {\sl contracted graph\/} $\Gamma/\gamma$ is the graph obtained by replacing all connected components of $\gamma$ by their residues inside $\Gamma$. As an example the residue of a graph $\Gamma$ is equal to $\Gamma/\Gamma$.
\dessin{55mm}{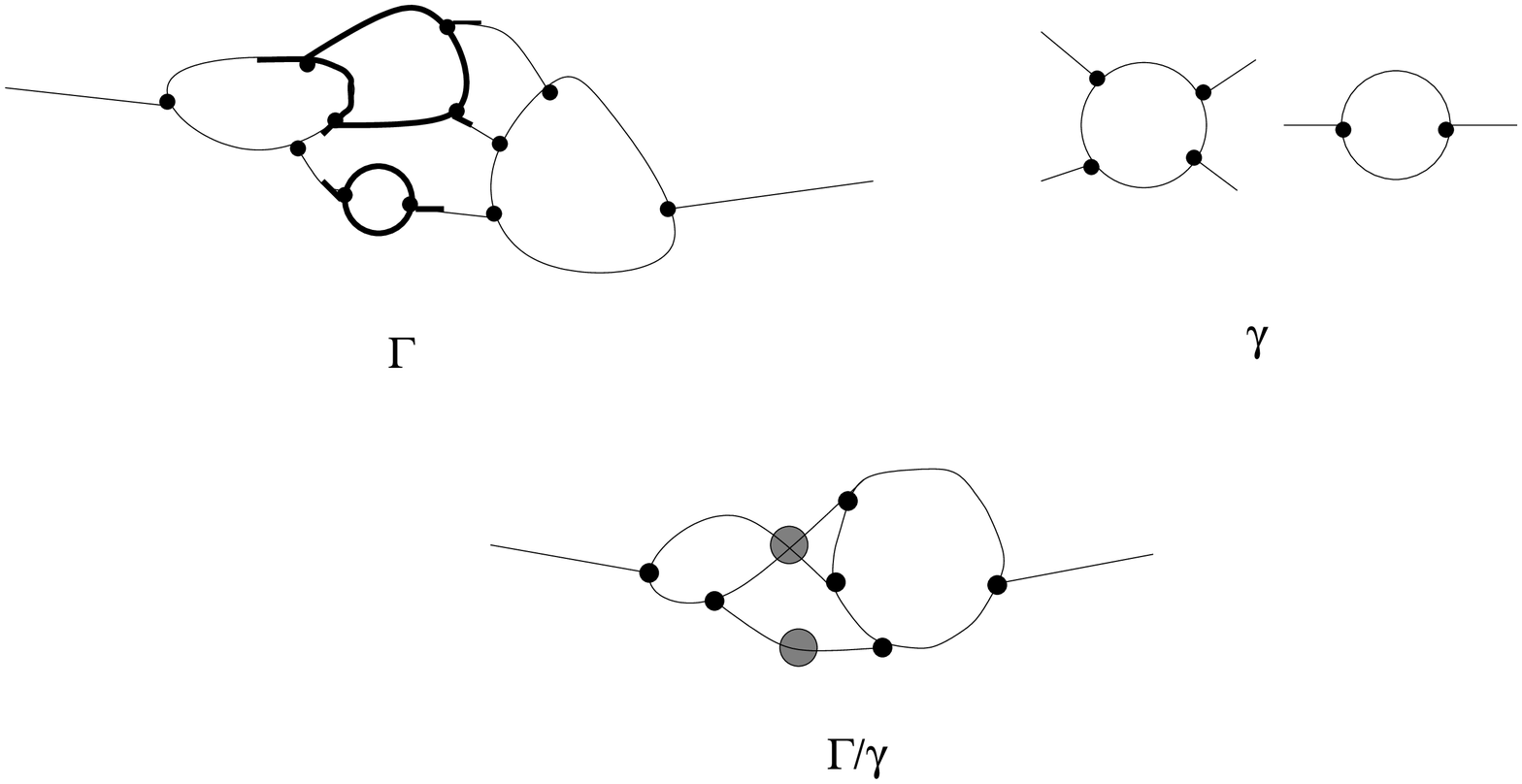}{4}{A subgraph $\scriptstyle \gamma$ inside a graph $\scriptstyle \Gamma$ in $\scriptstyle\varphi^{\scriptscriptstyle 3}$ theory. The contracted graph $\scriptstyle\Gamma/\gamma$ does not belong to $\scriptstyle\varphi^{\scriptscriptstyle 3}$}
\alinea{III.3. The graded Hopf algebra structure}
\qquad
Fix a set $\Cal T=\{\uple Tk\}$ of finite sequences of positive integers, which will be the possible vertex types we want to deal with. Let $V_{\Cal T}$ be the vector space generated by all connected 1PI Feynman graphs with vertex types in $\Cal T$, and all residues of those. Let $\Cal B_{\Cal T}=S(V_{\Cal T})$ be the free commutative algebra generated by $V$. We shall identify the unit $\un$ with the empty graph and any element of $\Cal B_{\Cal T}$ with a linear combination of disconnected locally 1PI graphs. The algebra structure is obvious, the co-unity is given by $\varepsilon(\un)=1$ and $\varepsilon(\Gamma)=0$ for any nonempty graph $\Gamma$. The grading (at least, one possible grading) is given on connected graphs by the {\sl loop number\/}~:
$$L:=I-V+1,$$
where $I$ is the number of internal edges and $V$ is the number of vertices of a given graph. This grading is extended to non-connected graphs in such a way that it is compatible with the algebra structure. It is important to notice that any nonempty subgraph has a non-vanishing loop number. The coproduct is given on connected 1PI graphs by the following formula~:
$$\eqalign{
\Delta(\Gamma)	&=\sum_{{\scriptstyle\gamma\hbox{ \sevenrm subgraph of }\Gamma
		\atop \scriptstyle \Gamma/\gamma\in V_\Cal T}}
\gamma\otimes \Gamma/\gamma	\cr
&=\Gamma\otimes\mop{res}\Gamma+\un\otimes\Gamma
	+\sum_{{\scriptstyle\gamma\hbox{ \sevenrm proper subgraph of }\Gamma
		\atop \scriptstyle \Gamma/\gamma\in V_\Cal T}}
\gamma\otimes \Gamma/\gamma \hskip 12mm\hbox{ if } L(\Gamma)\ge 1,\cr
\Delta(\Gamma)	&=\Gamma\otimes\Gamma\hskip 12mm\hbox{ if } L(\Gamma)=0,\cr}$$
and extended to non-connected graphs by multiplicativity. We leave it to the reader as an easy exercice to show that the coproduct respects the loop number as well. Figure 5 below illustrates a coproduct computation in $\varphi^3$ theory. Two terms of the sum have been removed because the corresponding contracted graphs have a vertex the type of which is outside $\Cal T$ (here a pentavalent and an hexavalent vertex respectively), and then does not belong to $V_{\Cal T}$. On the other hand residues with any number of external edges are allowed.
\vskip -6mm
\dessin{50mm}{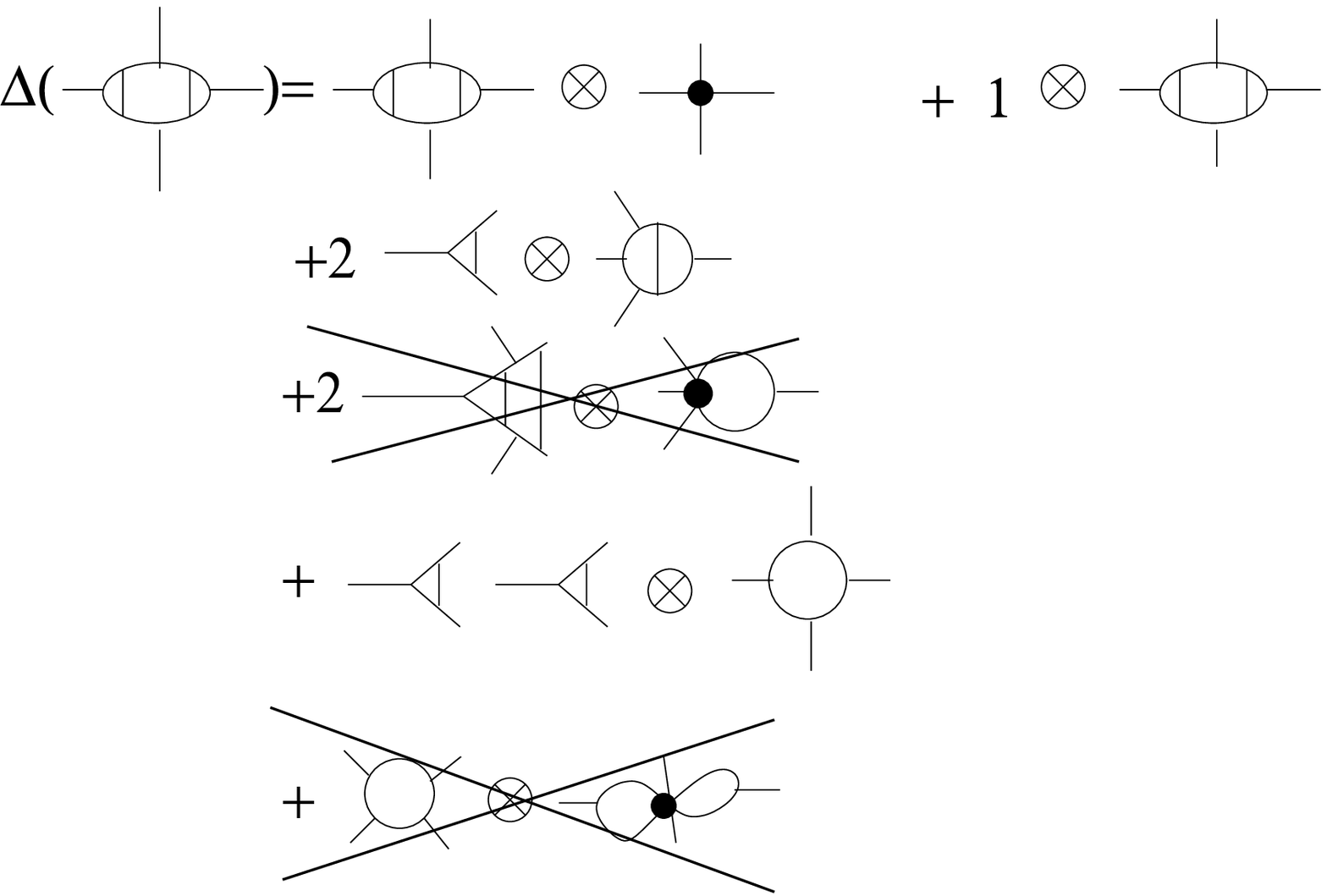}{5}{an example of coproduct in $\scriptstyle\varphi^{\scriptscriptstyle 3}$ theory}    
Figure 6 below illustrates another coproduct computation in $\varphi^3$ theory, with a bivalent vertex arising in the contracted graph~:
\dessin{18mm}{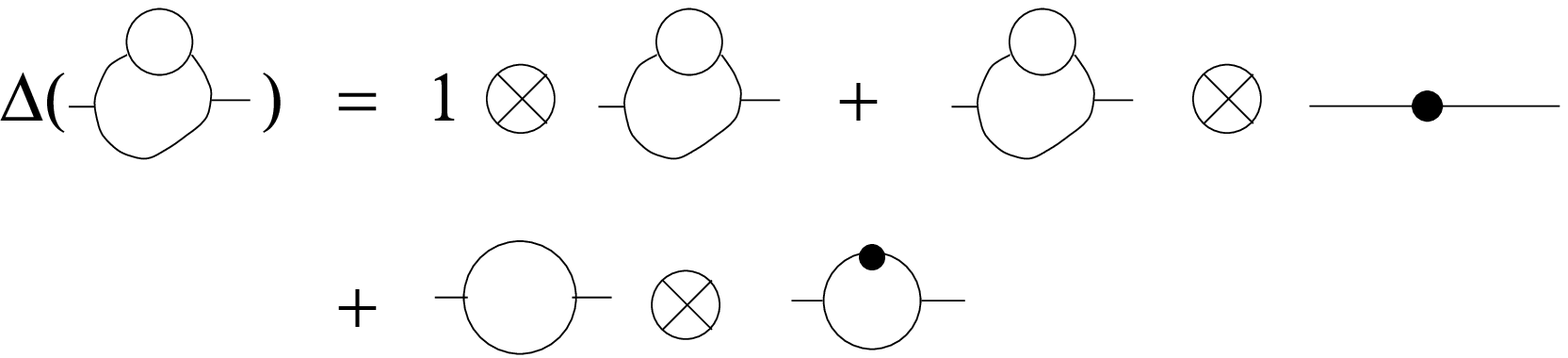}{6}{another example of coproduct in $\scriptstyle\varphi^{\scriptscriptstyle 3}$ theory}
\prop{III.3.1}
$\Cal B_{\Cal T}$ is a pointed graded bialgebra.
\dem
All axioms of a pointed graded bialgebra have been already given by the construction, except coassociativity of the coproduct. But we have for any 1PI graph of positive degree~:
$$(\Delta\otimes I)\Delta(\Gamma)=
\sum_{{\scriptstyle\delta\subset\gamma\subset\Gamma
		\atop \scriptstyle \gamma/\delta\in\Cal B_\Cal T,\ \Gamma/\gamma\in\Cal B_\Cal T}} \delta\otimes\gamma/\delta\otimes\Gamma/\gamma,$$
whereas~:
$$(I\otimes\Delta)\Delta(\Gamma)=
\sum_{{\scriptstyle \delta\subset\Gamma, \wt\gamma\subset\Gamma/\delta
\atop \scriptstyle \Gamma/\delta\in\Cal B_\Cal T,\ (\Gamma/\delta)/\wt\gamma\in\Cal B_\Cal T}} \delta\otimes\wt\gamma\otimes(\Gamma/\delta)/\wt\gamma.$$
There is an obvious bijection $\gamma\mapsto\wt\gamma=\gamma/\delta$ from subgraphs of $\Gamma$ containing $\delta$ onto subgraphs of $\Gamma/\delta$, given by shrinking $\delta$. As we have the obvious ``transitive shrinking property''~:
$$\Gamma/\gamma=(\Gamma/\delta)/\wt\gamma,$$
the two expressions coincide.
\qed
In order to build up a graded Hopf algebra from $\Cal B_{\Cal T}$, two choices are possible~: first we can add formally the inverses of the grouplike elements, i.e. the degree zero graphs~: let $\Sigma$ be the set of degree zero connected 1PI graphs, let $\Sigma\inver$ be another copy of the same set, with elements labelled $\gamma\inver,\gamma\in\Sigma$. Let $\wt V_{\Cal T}$ be the vector space generated by $V_{\Cal T}$ and $\Sigma\inver$, and consider~:
$$\wt\Cal H_{\Cal T}=S(\wt V_{\Cal T})/J,$$
where $J$ is the ideal generated by $\gamma\gamma\inver-\un,\,\gamma\in\Sigma$. The coproduct on $S(V_{\Cal T})$ is extended to $S(\wt V_{\Cal T})$ by saying that the elements of $\Sigma\inver$ are grouplike. The ideal $J$ is the a bi-ideal, and so $ \wt\Cal H_{\Cal T}$ is a pointed graded bialgebra. An antipode is easily given inductively with respect to the degree, as any degree zero element has an antipode given by $S(\gamma)=\gamma\inver$, $S(\gamma\inver)=\gamma$ for any $\gamma\in\Sigma$.
\ssq
The second option consists in killing the degree zero graphs (except the empty graph). We set~:
$$\Cal H_{\Cal T}=\Cal B_{\Cal T}/K,$$
where $K$ is the ideal generated by $\gamma-\un,\gamma\in\Sigma$. It is easily seen to be a bi-ideal. The quotient is then a connected graded bialgebra, hence a Hopf algebra thanks to corollary II.3.2. We can identify the quotient with $S(V'_{T})$, where $V'_T$ stands for the vector space generated by connected $1PI$ graphs with loop number $\ge 1$. The coproduct is then given by Kreimer's formula~:
$$\Delta(\Gamma)=\Gamma\otimes\un+\un\otimes\Gamma
	+\sum_{{\scriptstyle\gamma\hbox{ \sevenrm proper subgraph of }\Gamma
		\atop \scriptstyle \Gamma/\gamma\in V_\Cal T}}
\gamma\otimes \Gamma/\gamma.$$
\alinea{III.4. External structures\/}
\qquad We shall be very sketchy here. Let $W$ be a finite-dimensional vector space (the {\sl momentum space\/}). Keeping the notations of \S\ III.3 and following \cite{CK 2} and \cite {K 1}, a {\sl specified graph\/} will be a couple $(\Gamma,\sigma)$ where $\Gamma$ is a connected graph in $V_{\Cal T}$ with $E$ external lines, and $\sigma$ is a distribution on the vector subspace $M_\Gamma=M_E\subset W^E$ defined by~:
$$M_E=\{(\uple pE),\,\sum_{k=1}^Ep_k=0\}.$$
In order to get a Hopf algebra structure for specified graphs we must
discriminate further the type of a vertex~: once the number of edges of each
type is fixed for a vertex, we add an extra nonzero natural number, so that
there are ``several kinds of vertices of the same type''. This comes from the
lagrangian of the given quantum field theory we are dealing with~: each
monomial of degree $n_i$ with respect to the field $\phi_i$ $(i\in\{1,\ldots
,k\})$ gives rise to vertices of type $T=(\uple nk)$, and there are as many kinds of vertices of type $T$ as monomials of ``field degree'' $T$ inside the lagrangian. For example, in $\varphi^3$ theory with mass, terms $(m^2/2)\varphi^2$ and $(\partial \varphi)^2/2$ give rise to two different kinds of bivalent vertices. When taking residues we must specify the kind for the unique remaining vertex~: then when considering a contracted graph $\Gamma/\gamma$ we must consider the kind of every contracted vertex (corresponding to a connected component of the subgraph $\gamma$). This gives rise to contracted graphs $\Gamma/\gamma(i)$ where $i$ is a multi-index.
\ssq
To any vertex of kind $(T,i)$ corresponds a specific distribution $\sigma_{T,i}$ on $M_\Gamma$, where $\Gamma$ is any graph whose residue gives a vertex of type $T$. This extends to non-connected graphs by considering multi-indices $i$. Now $\Cal T$ stands for the set of all {\sl kinds\/} $(T,i)$ of vertices we can encounter, $V_{\Cal T}$  stands for the space generated by all connected 1PI graphs with vertex kinds in $\Cal T$, and all residues of those. Let $V'_{\Cal T}$ the space generated by all connected 1PI Feynman graphs with vertex kinds in $\Cal T$ and nonzero loop number, let $(V'_{\Cal T})_E$ the subspace of $V'_{\Cal T}$ of graphs with $E$ external edges, and finally let $W'_\Cal T$ the corresponding space of specified graphs~:
$$W'_\Cal T=\sum_{E=0}^\infty (V'_\Cal T)_E\otimes\Cal D'(M_E).$$
We directly give the connected version of the Hopf algebra~: it is given by $\Cal H_\Cal T=S(W'_{\Cal T})$, and the coproduct is given on connected specified graphs by~:
$$\Delta(\Gamma,\sigma)=(\Gamma,\sigma)\otimes\un+\un\otimes(\Gamma,\sigma)
	+\sum_{\scriptstyle\gamma\hbox{ \sevenrm proper subgraph of }\Gamma}
\ \sum_{i,\ \Gamma/\gamma(i)\in V_\Cal T}
(\gamma,\sigma_{T,i})\otimes (\Gamma/\gamma(i),\sigma).$$

\paragraphe{IV. An approach to the renormalization group}
\qquad Keeping the notations of paragraph II, we denote by $\Cal H$ a connected graded Hopf algebra, and by $\Cal A$ the algebra of germs of meromorphic functions at some $z_0\in\C$. The algebra $\Cal A$ admits a splitting into two subalgebras~:
$$\Cal A=\Cal A_+\oplus\Cal A_-,$$
where $\Cal A_+$ is the algebra of germs of holomorphic functions at $z_0$, and $\Cal A_-=(z-z_0)\inver\C[(z-z_0)\inver]$. We denote by $Y$ (resp. $\theta_t$) the biderivation (resp. the one-parameter group of automorphisms) of the Hopf algebra $\Cal H$ induced by the graduation (cf. \S\ II.6). As in paragraph II we denote by $G$ the group of the elements $\varphi\in\Cal L(\Cal H,\Cal A)$ such that $\varphi(\un)=\un_\Cal A$ (with the convolution products), and by $\g g$ the subalgebra of $\Cal L(\Cal H,\Cal A)$ of the elements $\varphi\in\Cal L(\Cal H,\Cal A)$ such that $\varphi(\un)=0$.
\ssq
Recall that $G=\exp \g g$. As in \S\ II we shall consider the subgroups $G_1$ (resp. $G_2$) of $G$ formed by the characters of $\Cal H$ with values in $\Cal A$ (resp. by the elements of $G$ which enjoy the cocycle property), as well as the Lie subalgebras $\g g_1$ (resp. $\g g_2$) of derivations of $\Cal H$ with values in $\Cal A$ (resp. of $\g g$ which enjoy the cocycle property). We have $G_1=\exp \g g_1$ and $G_2=\exp \g g_2$.
\alinea{IV.1. The renormalization map}	
We settle here a bijection $R:\g g\to\g g$ thanks to the biderivation $Y$~:
\prop{IV.1.1}
The equation~:
$$\varphi\circ Y=\varphi*\gamma \eqno{(E)}$$
defines a bijective correspondence~:
$$\eqalign{\wt R:G	&\longrightarrow \g g	\cr
		\varphi	&\longmapsto \gamma.	\cr}$$
Equivalently the equation~:
$$e^{*\alpha}\circ Y=e^{*\alpha}*\gamma \eqno{(E')}$$
defines a (non-linear) bijective correspondence~:
$$\eqalign{R:\g g	&\longrightarrow \g g	\cr
		\alpha	&\longmapsto \gamma,	\cr}$$
and $R=\wt R\circ\exp$.
\dem
Equation $(E)$ yields for any homogeneous $x\in\Cal H$~:
$$|x|\varphi(x)=\gamma(x)+\sum_{(x)}\varphi(x')\gamma(x''),$$
which determines $\gamma$ (recursively in $|x|$) from $\varphi$ and vice-versa, starting from $\varphi(\un)=\un_\Cal A$ and $\gamma(\un)=0$. In other words equation $(E)$ determines a bijection $\wt R$ from $G$ to $\g g$ such that $\gamma=\wt R(\varphi)$. The remainder of prop. IV.1.1 follows then immediately.
\qed
Equation $(E')$ yields the following explicit expression for $R$~:
$$R(\alpha)=e^{*-\alpha}*(e^{*\alpha}\circ Y).$$
There is another explicit formula~:
\prop{IV.1.2}
$$R(\alpha)=\int_0^1 e^{*-s\alpha}*(\alpha\circ Y)*e^{*s\alpha} \,ds.$$
\dem
for any $u\in\R$ we have~:
$$e^{*u\alpha}\circ Y=e^{*u\alpha}*R(u\alpha).$$
Setting $u=t+s$ and using the group property $e^{*(t+s)}\alpha=e^{*t\alpha}*e^{*s\alpha}$ as well as the derivation property~:
$$(e^{*t\alpha}*e^{*s\alpha})\circ Y=(e^{*t\alpha}\circ Y)*e^{*s\alpha}
		+e^{*t\alpha}*(e^{*s\alpha}\circ Y),$$
we get~:
$$e^{*(t+s)\alpha}\circ Y=e^{*(t+s)\alpha}*\bigl(R(s\alpha)
		+e^{*-s\alpha}*R(t\alpha)*e^{*s\alpha}\bigr).$$
Set $\gamma(t)=R(t\alpha)$~: the above equation reads~:
$$\gamma(t+s)=\gamma(s)+e^{*-s\alpha}*\gamma(t)*e^{*s\alpha}.$$
We have $\gamma(0)=0$, and differentiating this equation with respect to $s$ at $s=0$ yields~:
$$\dot\gamma(t)=\dot\gamma(0)+[\gamma(t),\alpha].$$
Differentiating once again with respect to $t$ gives then~:
$$\ddot \gamma(t)=[\dot\gamma(t),\alpha].$$
The solution of this first order differential equation is given by~:
$$\dot\gamma(t)=e^{*-t\alpha}*\dot\gamma(0)*e^{*t\alpha}.$$
Expanding the equation $e^{*t\alpha}\circ Y=e^{* t\alpha}*\gamma(t)$ up to order $1$ in $t=0$ yields immediately~:
$$\dot\gamma(0)=\alpha\circ Y.$$
Integrating and setting $t=1$ establishes then proposition IV.1.2.
\qed
\cor{IV.1.3}
Correspondence $R$ sends $\Cal A$-valued derivations to $\Cal A$-valued derivations and cocycles to cocycles.
\dem
First assertion follows immediately from propositions IV.1.2, II.4.2 and II.8.2. Second assertion follows directly from proposition IV.1.2.
\qed
{\sl Remark\/}~: If the Hopf algebra $\Cal H$ is cocommutative, then thanks to the commutativity of $\Cal A$, the convolution product is commutative. The correspondence $R$ becomes then linear and we simply have~:
$$R(\alpha)=\alpha\circ Y.$$
\alinea{IV.2. Inverting $\wt R$~: the scattering map\/}
We shall give an explicit expression of the map $\wt R\inver~:\g g\to G$. It takes the form~:
$$\wt R\inver(\gamma)=\mopl{lim}_{t\to +\infty} \exp{-tA}\exp{tB},$$
(cf. theorem IV.2.1 below), where $A$ and $B$ live in a semi-direct product
Lie algebra $\wt {\g g}=\g g\semi\C$. We have to describe this semi-direct product and the corresponding semi-direct product group $\wt G=G\semi\C$, and then we must endow $\wt G$ with a topology so that the above limit makes sense. We adapt here the proof of Theorem 2 in \cite{CK2}. To be precise, we define the Lie algebra~:
$$\wt g=\g g\semi\C.Z_0,$$
where the action of $Z_0$ on $\g g$ is given by the derivation~:
$$Z_0(\gamma)=\gamma\circ Y$$
(see lemma II.7.1). The corresponding group is $\wt G=G\semi\C$, where the right action of $\C$ on $G$ is given by~:
$$\varphi.t=\varphi\circ\theta_t,$$
so that the product is given by $(\varphi,t)(\psi,s)=\bigl(\varphi*(\psi\circ\theta_t),t+s\bigr).$ We shall not dig out a Lie group structure for $\wt G$ here, but we shall define the exponential map $\exp : \wt{\g g}\to\wt G$. It should of course coincide with the exponential already defined on $G$, and should verify~:
$$\exp tZ_0=(e,t)$$
so that $\exp tZ_0$ indeed acts on $G$ by composition with $\theta_t=\exp tY$
on the right. We should be able in principle to define $\exp (tZ_0+\gamma)$ by
means of the Baker-Campbell-Hausdorff formula as long as convergence problems
can be handled here. We prefer, like in \cite {CK2}, to give an alternative definition based on Araki's expansion formula \cite {Ar}~:
$$\hskip -5pt\exp(tZ_0+\gamma)	=\sum_{n=0}^\infty
\int_{\sum_{j=0}^n u_j=1,\, u_j\ge 0}
\exp(u_0tZ_0)\gamma\exp(u_1tZ_0)\gamma\cdots\gamma\exp(u_ntZ_0)
\,du_1\cdots du_n.$$
Let us check that the sum above makes sense in our particular context~: setting $v_j=u_j+u_{j+1}+\cdots+ u_n$ we get~:
$$\eqalign{&\exp (-tZ_0)\exp(tZ_0+\gamma)=\exp (-tZ_0).\cr
&\hskip -12mm\exp (tZ_0).\sum_{n=0}^\infty
\int_{0\le v_n\le\cdots\le v_1\le 1}\hskip -31 pt	
\exp (-tv_1Z_0)\gamma\exp (tv_1Z_0)\cdots
\exp (-tv_nZ_0)\gamma\exp (tv_nZ_0)\,dv_1\cdots dv_n	\cr
&=\sum_{n=0}^\infty
\int_{0\le v_n\le\cdots\le v_1\le 1}\hskip -31 pt	
(\gamma\circ\theta_{-tv_1})*\cdots*
(\gamma\circ\theta_{-tv_n})\,dv_1\cdots dv_n.	\cr}$$
The sum here is well defined as a locally finite sum, as it ends up at $n=n_0$ when evaluated at any $x=\Cal H^{n_0}$. It remains to check that the exponential thus defined enjoys the one-parameter group property. Indeed, for any $s,t$ real we have~:
$$\eqalign{\exp t(Z_0+\gamma)\exp s(Z_0+\gamma)
&=e^{tZ_0}\Bigl(\sum_{p=0}^\infty
t^p\int_{0\le v_p\le\cdots\le v_1\le 1}\hskip -31 pt	
(\gamma\circ\theta_{-tv_1})*\cdots*
(\gamma\circ\theta_{-tv_p})\,dv_1\cdots dv_p\Bigr).	\cr
&\hskip -40mm e^{sZ_0}\Bigl(\sum_{p=0}^\infty
s^q\int_{0\le w_q\le\cdots\le w_1\le 1}\hskip -31 pt	
(\gamma\circ\theta_{-sw_1})*\cdots*
(\gamma\circ\theta_{-sw_q})\,dw_1\cdots dw_q\Bigr)	\cr
&=e^{(t+s)Z_0}\sum_{p,q=0}^\infty t^ps^q\int\!\!\int_
{0\le v_p\le\cdots\le v_1\le 1,\, 0\le w_q\le\cdots\le w_1\le 1}\cr
&\hskip -40mm (\gamma\circ\theta_{-s-tv_1})*\cdots*(\gamma\circ\theta_{-s-tv_p})
*(\gamma\circ\theta_{-sw_1})*\cdots*
(\gamma\circ\theta_{-sw_q})\,dv_1\cdots dv_pdw_1\cdots dw_q	\cr
&=e^{(t+s)Z_0}\sum_{p,q=0}^\infty\int\!\!\int_
{0\le v_p\le\cdots\le v_1\le t,\, 0\le w_q\le\cdots\le w_1\le s}\cr
&\hskip -40mm (\gamma\circ\theta_{-s-v_1})*\cdots*(\gamma\circ\theta_{-s-v_p})
*(\gamma\circ\theta_{-w_1})*\cdots*
(\gamma\circ\theta_{-w_q})\,dv_1\cdots dv_pdw_1\cdots dw_q	\cr
&=e^{(t+s)Z_0}\sum_{n=0}^\infty\sum_{p+q=n}\int\!\!\int_
{s\le v_p\le\cdots\le v_1\le t+s,\, 0\le w_q\le\cdots\le w_1\le s}\cr
&\hskip -40mm (\gamma\circ\theta_{-v_1})*\cdots*(\gamma\circ\theta_{-v_p})
*(\gamma\circ\theta_{-w_1})*\cdots*
(\gamma\circ\theta_{-w_q})\,dv_1\cdots dv_pdw_1\cdots dw_q	\cr
&=e^{(t+s)Z_0}\sum_{n=0}^\infty\int_
{s\le u_p\le\cdots\le u_1\le t+s}\cr
&\hskip -40mm (\gamma\circ\theta_{-u_1})*\cdots*(\gamma\circ\theta_{-u_n})
\,du_1\cdots du_n	\cr
&=\exp (t+s)(Z_0+\gamma).
}$$
We can now state the main theorem of this section~:
\th{IV.2.1}
Let $\gamma\in\g g$. Then~:
\medskip
1) For any real $t$ the product $\exp -tZ_0\exp t(Z_0+\gamma)$ belongs to $G$.
\smallskip
2) The product above admits a limit when $t\to+\infty$ for the topology on $G$ induced by the simple convergence topology on $\Cal L(\Cal H,\Cal A)$.
\smallskip
3) The inverse of the renormalization map is given by~:
$$\wt R\inver(\gamma)=\mopl{lim}_{t\to +\infty}\exp -tZ_0\exp t(Z_0+\gamma).$$
\smallskip
4) $\wt R\inver$ sends $\g g_1$ into $G_1$ and $\g g_2$ into $G_2$.
\dem
The first assertion comes directly from the expression~:
$$\exp (-tZ_0)\exp(tZ_0+t\gamma)=\sum_{n=0}^\infty
\int_{0\le v_n\le\cdots\le v_1\le 1}\hskip -31 pt	
(t\gamma\circ\theta_{-tv_1})*\cdots*
(t\gamma\circ\theta_{-tv_n})\,dv_1\cdots dv_n.$$
The right-hand side belongs manifestly to $G$. Change of variables $v_j\to tv_j$ yields~:
$$\exp (-tZ_0)\exp(tZ_0+t\gamma)=\sum_{n=0}^\infty
\int_{0\le v_n\le\cdots\le v_1\le t}\hskip -31 pt	
(\gamma\circ\theta_{-v_1})*\cdots*
(\gamma\circ\theta_{-v_n})\,dv_1\cdots dv_n.$$
To prove the second assertion it suffices to prove that the integrals~:
$$I_n:=\int_{0\le v_n\le\cdots\le v_1\le +\infty}\hskip -31 pt	
(\gamma\circ\theta_{-v_1})*\cdots*
(\gamma\circ\theta_{-v_n})\,dv_1\cdots dv_n$$
converge, as the sum $I_0+I_1+I_2+\cdots$ is locally finite. The convergence is easily seen by induction on $n$~: indeed we have $I_0=e$ and the crucial equality valid for any $x\in\mop{Ker}\varepsilon$~:
$$Y\inver(x)=\int_0^\infty  \theta_{-t}(x)\,dt.$$
It follows that we have for any $a\in\g g$~:
$$\int_0^\infty a\circ\theta_{-t}\,dt=a\circ Y\inver.$$
A simple computation then gives~:
$$\eqalign{I_n	&=\int_0^\infty (I_{n-1}*\gamma)\circ\theta_{-v_n}\,dv_n	\cr
		&=(I_{n-1}*\gamma)\circ Y\inver,\cr}$$
which inductively shows the convergence of the integrals $I_n$. Now equation $(E)$ can be rewritten as~:
$$\eqalign{\varphi(x)	&=(\varphi*\gamma)\circ Y\inver(x)\ \ \forall x\in
\mop{Ker}\varepsilon	\cr
	\varphi(\un)	&=\un_{\Cal A}.	\cr}\eqno{(E'')}$$
As $\gamma=\wt R(\varphi)$ it means that~:
$$\wt R\inver (\gamma)=e+T\bigl(\wt R\inver(\gamma)\bigr),$$
where $T$ is the transformation of $\Cal L=\Cal L(\Cal H, \Cal A)$ defined by~:
$$\eqalign{T(\psi)	&=(\psi*\gamma)\circ Y\inver	\cr
			&=\int_0^\infty (\psi*\gamma)\circ\theta_{-t}\, dt.}$$
Transformation $T$ is a contraction on $\Cal L$ for the distance associated with the filtration. $\wt R\inver(\gamma)$ is then the limit of the sequence $(\varphi_n)$ defined by $\varphi_0=e$ and $\varphi_{n+1}=e+T(\varphi_n)$. A straightforward computation yields~:
$$\varphi_n=\sum_{k=0}^n I_k.$$
Hence we have~:
$$\wt R\inver(\gamma)=\sum_{k=0}^\infty I_k,$$
which proves assertion 3). Finally assertion 4) comes from the fact that derivation $Z_0$ acts on $\g g_1$ and $\g g_2$. We can then consider semi-direct products~:
$$\eqalign{&\wt{\g g}_1=\g g_1\semi\C.Z_0,\hskip 12mm
		\wt G_1=G_1\semi\C,	\cr
	&\wt{\g g}_2=\g g_2\semi\C.Z_0,\hskip 12mm
		\wt G_2=G_2\semi\C,	\cr}$$
and thus replace the group $G$ by any of the two groups $G_1,G_2$ in assertions 1), 2) and 3), which proves assertion 4) and ends the proof of the theorem. 
\qed
\cor{IV.2.2}
The inverse of $R:\g g\to\g g$ is given by~:
$$R\inver(\gamma)=\mopl{lim}_{t\to +\infty}\mop{Log}\bigl(
\exp -tZ_0\exp t(Z_0+\gamma)\bigr),$$
and $R\inver$ sends $\g g_1$ (resp. $\g g_2$) into $\g g_1$ (resp. $\g g_2$).
\ndem
\alinea{IV.3. The residue}
We keep the same notations as before except that we set $z_0=0$ for notational simplicity. To any $\psi\in\Cal L$ we associate a linear form $\mop{Res}\psi$ on $\Cal H$ by extracting the $z\inver$ term~: more precisely if we have for any $x\in\Cal H$ and for any $z$ in some pointed neighbourhood of $0$~:
$$\psi(x)(z)=\sum_{n=-N}^{+\infty}\psi_n(x)z^n$$
with $\psi_n(x)\in\C$, then~:
$$\mop{Res}\psi(x):=\psi_{-1}(x).$$
\alinea{IV.4. Renormalization map and Birkhoff decomposition}
Following more closely A. Connes and D. Kreimer \cite  {CK2} we shall focus our attention on elements of $G$ the Birkhoff decomposition of which shares an invariance property with respect to the action of the graduation. More precisely we define first a new action of $\C$ on the group $G$ by~:
$$\psi^t(x)(z)=\psi (\theta_{tz}(x))(z)=e^{tz|x|}\psi(x)(z).$$
We shall consider the elements $\psi$ of $G$ such that in the Birkhoff decomposition~:
$$\psi^t=(\psi^t)_-^{*-1}*(\psi^t)_+$$
the polar part $(\psi^t)_-$ is independent of $t$. The motivation to look at this very specific property is that it is fulfilled by the characters of Hopf algebras of Feynman graphs obtained via Feynman rules and dimensional regularization. In physical terms this corresponds to the well-known fact that the counterterms do not depend on the choice of the arbitrary mass $\mu$ one must introduce in order to perform dimensional regularization. Let us say that $\psi \in G$ enjoys property $(\Phi)$ if the above condition on the polar part of its Birkhoff decomposition is fulfilled. $G^\Phi$ will denote the subset of elements of $G$ enjoying property $(\Phi)$.
\ssq
For any $\psi\in\Cal L$ and any $f\in\Cal A$, we shall denote by $f\psi$ the element of $\Cal L$ defined by~:
$$\forall x\in\Cal H,\ f\psi(x)(z)=f(z)\psi(x)(z).$$
In particular, $(z\psi)$ is defined by $(z\psi)(x)(z)=z.\psi(x)(z)$.
\th{IV.4.1}
The map~:
$$\eqalign{z\wt R:G	&\longrightarrow	\g g	\cr
	\psi	&\longmapsto z\wt R(\psi)	\cr}$$
restricts to a bijection from $G^\Phi$ onto $\g g\cap\Cal L(\Cal H,\Cal A_+)$.
\dem
For any $\beta\in\g g$ introduce the linear transformation $U_\beta$ of $\g g$ defined by~:
$$U_\beta(A)=\beta*A+zA\circ Y.$$
If $\beta$ belongs to $\g g\cap\Cal L(\Cal H,\Cal A_+)$ then $U_\beta$ restricts to a linear transformation of $\g g\cap\Cal L(\Cal H,\Cal A_+)$.
\lemme{IV.4.2}
For any $\psi \in G, n\in\N$ we have~:
$$z^n\psi\circ Y^{n}=\psi*U_{z\wt R(\psi)}^n(e).$$
\dem
Case $n=0$ is obvious, $n=1$ is just the definition of $\wt R$. We check thus by induction~:
$$\eqalign{z^{n+1}\psi\circ Y^{n+1}	&=z(\psi\circ Y^n)\circ Y	\cr
				&=z\bigl(\psi*U_{z\wt R(\psi)}^n(e)\bigr)
\circ Y	\cr
				&=z(\psi\circ Y)*U_{z\wt R(\psi)}^n(e)
+z\psi*\bigl(U_{z\wt R(\psi)}^n(e)\circ Y\bigr)	\cr
				&=\psi*\bigl(z\wt R(\psi)*U_{z\wt R(\psi)}^n(e)
+zU_{z\wt R(\psi)}^n(e)\circ Y\bigr)	\cr
				&=\psi*U_{z\wt R(\psi)}^{n+1}(e).\cr}$$
\qed
Let us finish the proof of Theorem IV.4.1~: according to Lemma IV.4.2 we have for any real $t$, at least formally~:
$$\psi^t=\psi*\exp (tU_{z\wt R(\psi)})(e).\eqno{(*)}$$
We have still to fix the convergence of the exponential just above in the case when $z\wt R(\psi)$ belongs to $L(\Cal H,\Cal A_+)$. Let us consider the following decreasing bifiltration of $\Cal L(\Cal H,\Cal A_+)$~:
$$\Cal L_+^{p,q}=\bigl(z^q\Cal L(\Cal H,\Cal A_+)\bigr)\cap \Cal L^p.$$
Considering the associated filtration~:
$$\Cal L_+^n=\bigcup_{p+q=n}\Cal L_+^{p,q},$$
we see that for any $\beta\in\g g\cap\Cal L(\Cal H,\Cal A_+)$ transformation $U_\beta$ increases filtration by $1$, i.e~:
$$U_\beta(\Cal L^n_+)\subset\Cal L^{n+1}_+.$$
The algebra $\Cal L(\Cal H,\Cal A_+)$ is {\sl not\/} complete with respect to
the topology induced by this filtration, but the completion is  $\Cal L(\Cal H,\widehat{\Cal
  A_+})$, where $\widehat{\Cal
  A_+}=\C[[z]]$ stands for the  formal series. Hence the right-hand side of
(*) is convergent in $\Cal L(\Cal H,\widehat{\Cal
  A_+})$ with respect to
this topology. Hence for any $\gamma\in\Cal L(\Cal H,\Cal A_+)$ and for $\psi$
such that $z\wt R(\psi)=\gamma$ we have $\psi^t=\psi*h_t$ with $h_t\in\Cal
L(\Cal H,\widehat{\Cal A_+})$ for any $t$. On another hand we already know
that $h_t$ takes values in meromorphic functions for each $t$. So $h_t$
belongs to $\Cal L(\Cal H,\Cal A_+)$.
\ssq
Then for any $\beta\in\Cal L(\Cal H,\Cal A_+)$ and for $\psi$ such that $z\wt R(\psi)=\beta$ we have $\psi^t=\psi*h(t)$ with $h(t)\in\Cal L(\Cal H,\Cal A_+)$ for any $t$. This is equivalent to the fact that $\psi$ belongs to $G^\Phi$. Conversely take any $\psi$ in $G^\Phi$. Then there exists $h_t\in\Cal L(\Cal H,\Cal A_+)$ such that $\psi^t=\psi*h_t$. Hence~:
$${d\over dt}\restr{t=0}\psi_t=z(\psi\circ Y)=\psi\dot h_t\restr{t=0}.$$
So $z\wt R(\psi)=\dot h_t\restr{t=0}$ belongs to $\Cal L(\Cal H,\Cal A_+)$. This proves theorem IV.4.1.
\qed 
\lemme{IV.4.3}
1) Let $\psi\in G^\Phi$, and let~:
$$\psi=(\psi_-)^{*-1}*\psi_+$$
its Birkhoff decomposition. Then $\wt\psi=(\psi_-)^{*-1}$ enjoys property $(\Phi)$ as well.
\smallskip
2) Let $\psi\in G^\Phi$ and let $h$ be any element of $G$ without polar part. Then $\psi*h$ enjoys property $(\Phi)$ as well.
\dem
Assertion 1) is a consequence of the equality~:
$$\psi^t=(\psi^t)_-^{*-1}*(\psi^t)_+=\wt\psi^t*(\psi_+)^t,$$
which implies immediately~:
$$\wt\psi^t=(\psi^t)_-^{*-1}*(\psi^t)_+*
	\bigl((\psi_+)^t\bigr)^{*-1}.$$
This is the Birkhoff decomposition of $\wt\psi^t$~: its polar part is $(\psi^t)_-$, which by hypothesis is independant of $t$. Assertion 2) comes from the equality~:
$$(\psi*h)^t=\psi^t*h^t=(\psi^t)_-^{*-1}*\bigl((\psi^t)_+*h^t\bigr),$$
and from the fact that $h_t$ does not have any polar part. Thus the right-hand side is the Birkhoff decomposition of $(\psi*h)^t$, the polar part of which is clearly independant of $t$. 
\qed
Thanks to Lemma IV.4.3 we can now focus our attention on elements $\psi$ of $G$ which enjoy property $(\Phi)$ and such that $\psi$ sends $\mop{Ker}\varepsilon$ into $\Cal A_-$. Let us denote by $G^\Phi_-$ the subset of $G$ formed by these elements.
\th{IV.4.4}
Let $\g g^c=\{\beta\in\Cal L(\Cal H,\C),\ \beta(\un)=0\}$. Let $\g g^c_1=\g g^c\cap\g g_1$, and $\g g^c_2=\g g^c\cap\g g_2$. Let $G^\Phi_{1,-}= G^\Phi_-\cap G_1$ and $G^\Phi_{2,-}= G^\Phi_-\cap G_2$. Then the map~:
$$\eqalign{G	&\longrightarrow	\g g	\cr
	\psi	&\longmapsto z\wt R(\psi)	\cr}$$	
restricts to a bijection from $G^\Phi_-$ (resp. $G^\Phi_{1,-}$, $G^\Phi_{2,-}$) to $\g g^c$ (resp $\g g^c_1$, $\g g^c_2$) , and we have explicitly for any $\psi\in G^\Phi_-$~:
$$\wt R(\psi)={1\over z}\bigl((\mop{Res}\psi)\circ Y\bigr).$$
\dem
We shall need the following key lemma~:
\lemme{IV.4.5}
For any $\psi\in G^\Phi_-$ we have~:
$${d\over dt}\restr{t=0}\bigl((\psi^t)_+)=\mop{Res}(\psi\circ Y).$$
\dem
Using property $(\Phi)$, the fact that $\psi(\mop{Ker}\varepsilon)\subset\Cal A_-$ and the explicit expression of $(\psi^t)_+$ given by theorem II.4.1  we have for any $x\in\Cal H$~:
$$\eqalign{(\psi^t)_+(x)&=(I-\pi)\Bigl(
\psi^t(x)+\sum_{(x)}\psi^{*-1}(x')\psi^t(x'')\Bigr)	\cr
			&=t(I-\pi)\bigl(
z|x|\psi(x)+z\sum_{(x)}\psi^{*-1}(x')\psi(x'')|x''|\Bigr)+O(t^2)\cr
			&=t\mop{Res}\bigl(\psi\circ Y)+O(t^2).\cr
}$$
\qed
Now for any $\psi\in G^\Phi_-$ the Birkhoff decomposition of $\psi^t$ reads~:
$$\psi^t=\psi*(\psi^t)_+.$$
Differentiating with respect to $t$ at $t=0$ we get according to Lemma IV.4.5~:
$$z\psi\circ Y=\psi*\mop{Res}(\psi\circ Y).$$
We deduce then~:
$$\psi\circ Y=\psi*{1\over z}\mop{Res}(\psi\circ Y),$$
which proves the equality $\wt R(\psi)={1\over z}\mop{Res}(\psi\circ Y)$. As a consequence correspondence $z\wt R$ sends $G^\Phi_-$ into $\g g^c$. Conversely let $\beta$ in $\g g^c$. Consider $\psi=\wt R\inver(z\inver\beta)$. This element of $G$ verifies by definition~:
$$z\psi\circ Y=\psi*\beta.$$
Hence for any $x\in\mop{Ker}\varepsilon$ we have~:
$$z\psi(x)={1\over |x|}\Bigl(\beta(x)+\sum_
{(x)}\psi(x')\beta(x'')\Bigr).$$
As $\beta(x)$ is a constant (as a function of the complex variable $z$) it is easily seen by induction on $|x|$ that the right-hand side evaluated at $z$ has a limit when $z$ tends to infinity. Thus $\psi(x)\in\Cal A_-$, and then~:
$$\psi =\wt R\inver({1\over z}\beta)\in G^\Phi_-.$$
The proof of the variants of theorem IV.4.4 with $G_1$ and $G_2$ is then immediate thanks to corollary IV.1.3 and theorem IV.2.1.
\qed
{\bf Remark}~: the ``$G_1$'' version of Theorem IV.4.4 recovers the result of
A. Connes and D. Kreimer in \cite {CK2 \S\ 2}~: when $\psi$ is a character,
then the $\C$-valued derivation $z\wt R(\psi)$ is denoted by $\beta$ in {\it
loc. cit\/}. Theorem IV.4.4 appears then as an interpretation of the result of
Connes and Kreimer in other groups than the group of $\Cal A$-valued characters.%
\paragraphe{References}
\bib{Ab}E. Abe, {\sl Hopf Algebras\/}, Cambridge Univ. Press (1980).
\bib{ABS}M. Aguiar, N. Bergeron, F. Sottile, {\sl Combinatorial Hopf algebras
and generalized Dehn-Sommerville relations\/}, Compos. Math.  142  no. 1, 1-30.(2006).
\bib{Ar}H. Araki, {\sl Expansional in Banach algebras\/}, Ann. Scient. Ec. Norm. Sup. 4e s\'erie, {\bf 6}, 67-84 (1973).
\bib{B}N. Bourbaki, {\sl Alg\`ebre, Chapitre 8\/}, Hermann, Paris.
\bib{BF}Ch. Brouder, A. Frabetti, {\sl Noncommutative renormalization for
massless QED\/}, hep-th/0011161 (2000).
\bib{BP}N.N. Bogoliubov, O.S. parasiuk, {\sl On the multiplication of causal
functions in the quantum theory of fields\/}, Acta Math. {\bf 97}, 227-266 (1957).
\bib{C}J. Collins, {\sl Renormalization\/}, Cambridge (1984).
\bib{CK1}A. Connes, D. Kreimer, {\sl Renormalization in Quantum Field Theory and the Riemann-Hilbert problem I : The Hopf algebra structure of graphs and the main theorem\/}, Commun. Math. Phys. 210, 249-273 (2000).
\bib{CK2}A. Connes, D. Kreimer, {\sl Renormalization in quantum field theory and the Riemann-Hilbert problem. II. The $\beta$-function, diffeomorphisms and the renormalization group\/}. Comm. Math. Phys. 216, no. 1, 215--241 (2001).
\bib{CK3}A. Connes, D. Kreimer, {\sl Insertion and elimination: the doubly
infinite Lie algebra of Feynman graphs\/}, Ann. Henri Poincaré 3, no. 3,
411--433 (2002).
\bib{CM1}A. Connes, M.Marcolli, {\sl From physics to number theory via
noncommutative geometry II\/},
arxiv:math.QA/0411114 (2004).
\bib{CM2}A. Connes, M.Marcolli, {\sl A walk in the noncommutative garden\/},
arxiv:math.QA/0601054 (2006).
\bib{Di}J. Dixmier, {\sl Alg\`ebres enveloppantes\/}, Gautier-Villars, Paris (1974).
\bib{DF}R.K. Dennis, B. Farb, {\sl Noncommutative algebra\/}, Springer Verlag (1993).
\bib{DK}Yu. A. Drozd, V.V. Kirichenko, {\sl Finite dimensional algebras (english edition)\/}, Springer Verlag (1994).
\bib{DNR}S. D\v asc\v alescu, C. N\v ast\v asescu, S. Raianu, {\sl Hopf algebras, an introduction\/}, Pure and Applied mathematics vol. 235, Marcel Dekker (2001).
\bib{E}P. Etingof, {\sl Note on dimensional regularization\/}, in {\sl Quantum
fields and strings : a course for mathematicians\/}, AMS/IAS (1999).
\bib{EGK1}K. Ebrahimi-Fard, L. Guo, {\sl Matrix representation of of
renormalization in perturbative  quantum field theory\/}, preprint
arxiv:hep-th/0508155 (2005).
\bib{EGK1}K. Ebrahimi-Fard, L. Guo, D. Kreimer, {\sl Integrable renormalization I: the ladder case\/}, arXiv:hep-th/0402095 (2004).
\bib{EGK2}K. Ebrahimi-Fard, L. Guo, D. Kreimer, {\sl Integrable renormalization II: the general case\/}, arXiv:hep-th/0403118 (2004).
\bib{EGK3}K. Ebrahimi-Fard, L. Guo, D. Kreimer, {\sl Spitzer's identity and
the algebraic Birkhoff decomposition in pQFT\/}, J. Phys. A: Math. Gen. {\bf
37}, 11036-11052 (2004), arXiv:hep-th/0407082.
\bib{EGM}K. Ebrahimi-Fard, L. Guo, D. Manchon, {\sl Birkhoff type
decompositions and the Baker-Campbell-Hausdorff recursion\/},
Comm. Math. Phys. (to appear). arxiv:math-ph/0602004 (2006).
\bib{EGGV}K. Ebrahimi-Fard, J.M. Gracia-Bondia, L. Guo, J.C. V\`arilly, {\sl
Combinatorics of renormalization as matrix calculus\/}, Phys. Lett. B, {\bf
632} No 4, 552-558 (2006), arxiv:hp-th/0508154.
\bib{F}L. Foissy, {\sl Les alg\`ebres de Hopf des arbres enracin\'es d\'ecor\'es I,II\/}, Bull. Sci. Math. 126, 193-239 et 249-288 (2002).
\bib{FG}H. Figueroa, J.M. Gracia-Bond\'\i a, {\sl Combinatorial Hopf algebras
in Quantum Field Theory I\/}, Reviews of Mathematical Physics {\bf 17},
881-976 (2005).
\bib{GZ}L. Guo, B. Zhang, {\sl Renormalization of multiple zeta values\/},
preprint (2006).
\bib{He}K. Hepp, {\sl Proof of the bogoliubov-Parasiuk theorem on
renormalization\/}, Comm. Math. Phys. {\bf 2}, 301-326 (1966).
\bib{H}M.E. Hoffman, {\sl The Hopf algebra structure of multiple harmonic
sums\/}, Nuclear Phys. B Proc. Suppl.  135, 215-219  (2004). arXiv:math.QA/04\-06589.
\bib{I1}L.M. Ionescu, {\sl Perturbative quantum field theory and configuration space integrals\/}, arXiv:\-hep-th/0307062 (2003). 
\bib{I2}L.M. Ionescu, {\sl A combinatorial approach to coefficients in deformation quantization\/}, arXiv:\-hep-th/0404389 (2004).
\bib{IM}L.M. Ionescu, M. Marsalli, {\sl A Hopf algebra deformation approach to renormalization\/}, arXiv:hep-th/0307112 (2003).
\bib{J}N. Jacobson, {\sl Basic algebra II (second edition)\/}, Freeman, New-York, 1989.
\bib{Ka}C. Kassel, {\sl Quantum groups\/}, Springer Verlag (1995).
\bib{K1}D. Kreimer, {\sl Structures in Feynman graphs- Hopf algebras and symmetries\/}, hep-th/0202100 (2002).
\bib{K2}D. Kreimer, {\sl New mathematical structures in renormalizable quantum field theories\/}, hep-th/0211136 (2002).
\bib{Ma}D. Manchon, {\sl L'alg\`ebre de Hopf bitensorielle\/}, Comm. Algebra
25, no. 5, 1537--1551 (1997).
\bib{MP1}D. Manchon, S. Paycha, {\sl Shuffle relations for regularised
integrals of symbols\/}, arXiv:math-ph/0510067 (2005).
\bib{MP2}D. Manchon, S. Paycha, {\sl Renormalised Chen integrals for symbols
on $\R^n$ and renormalised polyzeta functions\/}, arXiv:math-ph/0604562 (2006).
\bib{Mo}S. Montgomery, {\sl Some remarks on filtrations of Hopf algebras\/},
Comm. Algebra 21 No 3, 999-1007 (1993).
\bib{R}H. Ratsimbarison, {\sl Feynman diagrams, Hopf algebras and
renormalization\/}, arxiv:math-ph/0512012 (2005).
\bib{RV}M. Rosenbaum, J.D. Vergara, {\sl The Hopf algebra of renormalization,
normal coordinates and Kontsevich deformation quantization\/},
arXiv:hep-th/0404233 (2004).
\bib{S}M. Sakakibara, {\sl On the Differential equations of the characters for
the Renormalization group\/}, Mod.Phys.Lett. A19, 1453-1456 (2004).
\bib{Sw}M. E. Sweedler, {\sl Hopf algebras\/}, Benjamin, New-York (1969).
\bib{TW}E.J. Taft, R.L. Wilson, {\sl On antipodes in pointed Hopf algebras\/},
J. Algebra 28, 27-32 (1974).
\bib{Z}W. Zimmermann, {\sl Convergence of Bogoliubov's method of
renormalization in momentum space\/}, Comm. Math. Phys. {\bf 15}, 208-234 (1969).

\bye